\documentclass[12pt]{article}
\usepackage{amssymb}
\usepackage{amsmath,color}
\usepackage[colorlinks,linkcolor=blue]{hyperref}
\usepackage{mathtools}

\begin{document}
	\title{\bf Global classical solutions to the compressible micropolar
		viscous fluids with large oscillations and vacuum}
			\author{\small Canze Zhu$^{*}$ ~and~ Qiang Tao\\
		{\small College of Mathematics and Statistics, Shenzhen University, Shenzhen 518060, China}	\thanks{Corresponding author.~canzezhu@163.com (C. Zhu).}}	
	\date{}
	\maketitle
	\date{}
	\maketitle
	\newtheorem{theorem}{Theorem}[section]
	\newtheorem{definition}{Definition}[section]
	\newtheorem{lemma}{Lemma}[section]
	\newtheorem{proposition}{Proposition}[section]
	\newtheorem{corollary}{Corollary}[section]
	\newtheorem{remark}{Remark}[section]
	\renewcommand{\theequation}{\thesection.\arabic{equation}}
	\catcode`@=11 \@addtoreset{equation}{section} \catcode`@=12
	{\bf Abstract.}{\small In this paper, we consider the three dimensional Cauchy problem of the compressible micropolar viscous flows, we prove the existence of unique global classical solution for smooth initial data with small initial energy but possibly large oscillations, the initial density may allowed to contain vacuum states. Furthermore, the large-time behavior of the solution is obtained.}
	
	\section{Introduction}	In this paper, we investigate the Cauchy problem for the compressible compressible micropolar viscous system \cite{1}, that is,
		\begin{align}\label{1.1}
		\begin{cases}
		\rho_{t}+\mathrm{div}(\rho \mathbf{u})=0, \\
		(\rho \mathbf{u})_{t}+\mathrm{div}(\rho \mathbf{u}\otimes \mathbf{u})+\nabla P(\rho)=(\mu+\zeta)\Delta \mathrm{\mathbf{u}}+(\mu+\lambda-\zeta)\nabla \mathrm{div}\mathbf{u}+2\zeta\mathrm{curl}\mathbf{w}, \\
		(\rho \mathbf{w})_{t}+\mathrm{div}(\rho \mathbf{u}\otimes \mathbf{w})+4\zeta\mathbf{w}=\mu^{'}\Delta{\mathbf{w}}+(\mu^{'}+\lambda^{'})\nabla \mathrm{div}\mathbf{w}+2\zeta\mathrm{curl}\mathbf{u},  \\
		\end{cases}
		\end{align}
		with $x \in \mathbb{R}^3$ and $t\geq 0$, where $\rho, \mathrm{\mathbf{u}}$ and $\mathrm{\mathbf{w}}$
		 represent density, velocity and microrotational velocity, respectively. $P(\rho)= \rho^\gamma$ with $\gamma>1$ represents the pressure. The positive constant
		viscosity coefficients $\mu$, $\lambda$, $\zeta$, $\mu^{'}$ and $\lambda^{'}$ satisfy the physical restrictions:
		\begin{align*}
		2\mu+3\lambda-4\zeta \ge 0 \quad \text{and} \quad 2\mu^{'}+3\lambda^{'}\ge 0.
		\end{align*}
		We consider the Cauchy problem of the system $(\ref{1.1})$ with initial data	
		\begin{equation}\label{1.2}
		(\rho,\mathbf{u},\mathbf{w})(x,t)|_{t=0} = \big(\rho_{0}(x),\mathbf{u}_{0}(x),\mathbf{w}_{0}(x)\big)
		\end{equation}
		and the far field behavior
		\begin{equation}\label{1.3}
		\quad(\rho,\mathbf{u},\mathbf{w})(x,t)=(\tilde{\rho},0,0)  ~~\text{as~} |x|\to\infty,~ t\ge 0.
		\end{equation}
		
		Compressible micropolar viscous fluids describing the viscous
		compressible fluids with randomly oriented particles suspended in the medium when the deformation of the fluid particles is ignored.
		Due to the effect of microparticles, this system show some physical phenomena which cannot be treated by the compressible viscous baratropic flows. Thus, it has mechanically significant for the microstructure of the polar fluids, see e.g.\cite{1,2,3}.
		
		There has been much works about the compressible micropolar viscous fluids. For the case of one-dimensional space, it has many results on the local or global existence, uniqueness and large time behavior of the solution (see e.g.\cite{5}-\cite{7}). The global existence of strong solutions to the 1D model with initial vacuum was obtained in \cite{4}. \cite{8} established the global existence and large time behavior of solution to the 2D micropolar equations with only angular viscosity dissipation. For the 3D model, the compressible micropolar fluid in the 3D case has been considered in by Dražić and Mujaković in the spherically symmetric case. Chen \cite{11} showed the local existence and uniqueness of strong solutions under the assumption that the initial density with vacuum. Liu and Zhang \cite{12} obtained global existence under the condition that the initial perturbation small in $H^N$ $(N\ge4)$.
		The blow up criterion of strong solutions to the Cauchy problem and the global weak solutions with discontinuous initial data and vacuum were proved by Chen et.al \cite{13} and \cite{14}. 
		
		For the related models of (\ref{1.1}), the magneto-micropolar fluid,
		Yuan showed the blow-up criterion for the smooth solution and gave the local smooth solution
		without the smallness of the initial data in \cite{15}. Amirat and Hamdache \cite{16} considered the global weak solutions with
		finite energy and establish the long-time behavior of the solution. Recently, Wei.et.al \cite{17} established the global existence under the condition that the initial data small in $H^3$. 
		
		Note that $\mathbf{w}= 0$ and $\zeta= 0$, (\ref{1.1}) reduces to the celebrated Navier–Stokes equation, it has many results. For Local existence, Choe and Kim established local in time strong solution of isentropic compressible fluids under the condition that initial density may vanish in an open subset in \cite{24}, Cho and Kim \cite{25} studied the local existence of strong solution with vacuum. For global existence, It is worth mentioning that \cite{18} first established global existence and uniqueness of solution for compressible Navier-Stokes equations with vacuum and smooth initial data which is of small energy but possibly large oscillations. Further relevant results, we can see \cite{21,22,23} and the references therein. 
	
		In this paper, we consider an initial value problem of the compressible micropolar
		viscous flows  $(\ref{1.1})$--$(\ref{1.3})$. Motivated by the works for compressible Navier-Stokes equation
		\cite{18} and the compressible MHD equation \cite{19}, we will first established global
		existence and uniqueness of solution with smooth initial data which is of small energy but possibly large oscillations, and the initial density may allowed to contain vacuum states. Then, the large time behavior of the solution will be given as well.
		we denote
		\begin{align*}
			\int f dx =\int_{\mathbb{R}^3} f dx.
		\end{align*}
		For $1\le\gamma\le\infty$ and $\beta>0$, we denote the standard homogeneous and inhomogeneous Sobolev spaces as follows:
		\begin{align*}
		\begin{cases}
		& L^r=L^r(\mathbb{R}^3),\ D^{k,r}=\{\mathbf{u}\in L^{1}_{loc}(\mathbb{R}^3)|\ ||\nabla^k\mathbf{u}||_{L^r}<\infty\}, \ ||\mathbf{u}||_{D^{k,r}}\triangleq||\nabla^k\mathbf{u}||_{L^r},\\
		& W^{k,r}=L^r\cap D^{k,r},\ H^{k}=W{k,2},\ D^k=D^{k,2},\ D^1=\{\mathbf{u}\in L^6 |\ ||\nabla\mathbf{u}||_{L^2}<\infty\},\\
		& \dot{H}^{\beta}=\{f:\mathbb{R}^3\rightarrow\mathbb{R}|\ ||f||_{ \dot{H}^{\beta}}^2=\int|\xi|^{2\beta}|\hat{f}(\xi)|^2d\xi<\infty \}.
		\end{cases}
		\end{align*}
		The initial energy is defined as
		\begin{align*}
			E_0= \int\bigg(\frac{1}{2}\rho_0|\mathbf{u}_0|^2+G(\rho_0)\bigg)dx,
		\end{align*}
		where $G$ denotes the potential energy density given by
		\begin{align*}
			G(\rho)\triangleq\rho\int_{1}^{\rho}\frac{P(s)-P(1)}{s^2}ds,
		\end{align*}
		It is obvious that
		\begin{align*}
		\begin{cases}
		 G(\rho)=\frac{1}{\gamma-1}P(\rho)\quad &\text{if}\ \tilde{\rho}=0,\\
		 G(\rho)\le C(\bar{\rho})(\rho-1)^2\quad &\text{if}\ \tilde{\rho}>0,\ 0\le\rho\le 2\bar{\rho},
		\end{cases}		
		\end{align*}
		for positive constant $C(\bar{\rho})$ depending on $\bar{\rho}$ and $\tilde{\rho}$.
		
		The main results in this paper are stated as follows.
		\begin{theorem}\label{t1.1}
			For given $M>0$, $\bar{\rho}\ge\tilde{\rho}+1$ and $q\in(3,6)$, assume that the initial data $(\rho_0,\mathbf{u}_0,\mathbf{w}_0)$ satisfy
			\begin{align}\label{1.4}
				\begin{cases}
				& G(\rho_0)+\rho_0|\mathbf{u}_0|^2+\rho_0|\mathbf{u}_0|^2+\rho_0|\mathbf{w}_0|^2\in L_1,\quad 0\le\rho_0\le\bar{\rho},\\
				&(\rho_0-\tilde{\rho},P(\rho_0)-P(\tilde{\rho}))\in H^2\cap W^{2,q},\\
				&\mathbf{u}_0\in D^1\cap D^2,\quad \mathbf{w}_0\in D^1\cap D^2,\\
				&||\mathbf{u}_0||_{H^1}+||\mathbf{w}_0||_{H^1}\le M,
				\end{cases}
			\end{align}and that the compatibility condition holds
			\begin{align}\label{1.5}
			&-(\mu+\zeta)\Delta \mathbf{u}_0-(\mu+\lambda-\zeta)\nabla \mathrm{div}\mathbf{u}_0-2\zeta\mathrm{curl}\mathbf{w}_0+\nabla P(\rho_0)=\rho_0^\frac{1}{2}g_1, \\
            &-\mu^{'}\Delta{\mathbf{w}}-(\mu^{'}+\lambda^{'})\nabla \mathrm{div}\mathbf{w}-2\zeta\mathrm{curl}\mathbf{u}+4\zeta\mathbf{w}=\rho_0^\frac{1}{2}g_2,			
			\end{align}for some $g_i\in L^2\ (i=1,2)$. Then, the Cauchy problem $(\ref{1.1})$-$(\ref{1.3})$ has a unique global classical solution $(\rho,\mathbf{u},\mathbf{w})$ satisfying
			\begin{align}\label{1.6}
			0\le\rho \le 2\bar{\rho}\quad \text{for all} \quad x\in\mathbb{R}^3, t\ge 0 ,
			\end{align}
			and 
			\begin{align}\label{1.7}
				\begin{cases}
			&(\rho-\tilde{\rho},P(\rho)-P(\tilde{\rho})\in C([0,T];H^2\cap W^{2,q}),\\
			&\mathbf{u}\in C([0,T];D^1\cap D^2)\cap L^\infty ([\tau,T];D^2\cap D^{3,q}),\\
			&\mathbf{u}_t\in L^\infty([\tau,T];D^1\cap D^2)\cap H^1([\tau,T];D^1),\\
			&\mathbf{w}\in C([0,T];D^1\cap D^2)\cap L^\infty ([\tau,T];D^2\cap D^{3,q}),\\
			&\mathbf{w}_t\in L^\infty([\tau,T];D^1\cap D^2)\cap H^1([\tau,T];D^1),		
				\end{cases}
			\end{align}
			for any $0<\tau<T<\infty$, provided that 
			\begin{align}
			E_0\le \epsilon
			\end{align} with constant $\epsilon$ depending only on $\mu, \lambda, \zeta, \mu^{'}, \lambda^{'},
			\gamma, \bar{\rho}, \tilde{\rho}\ \text{and}\ M$. Furthermore, the following large-time behavior holds:
			\begin{align}\label{1.10}
				\lim\limits_{t\to\infty}\bigg(||\rho-\tilde{\rho}||_{L^p}+\int\rho^\frac{1}{2}|\mathbf{u}|^4dx+||\nabla\mathbf{u}||_{L^r}+||\nabla\mathbf{w}||_{L^r}\bigg)=0
			\end{align}for $r\in[2,6)$ and 
			\begin{align}\label{1.8}p\in
				\begin{cases}
				&(\gamma,\infty)\quad\text{if}\quad\tilde{\rho}=0,\\
				&(2,\infty)\quad\text{if}\quad\tilde{\rho}>0.
				\end{cases}
			\end{align}
		\end{theorem}
	The rest of this paper is organized as follow. In section 2, we first establish a priori estimates of  of the smooth solution
	for problem $(\ref{1.1})$--$(\ref{1.3})$ with small initial energy.
	In section 3, we prove the Theorem \ref{t1.1}.
		\section{Preliminaries}
		We give the Sobolev inequalities (\cite{26}) at first.
		\begin{lemma}\label{l2.1}
			For $p\in[2,6]$, $q\in(1,\infty)$, and $r\in(3,\infty)$, there exists some feneric constant $C>0$ that may depend on $q$ and $r$ such that for $f\in H^1({\mathbb{R}^3})$ and $g\in L^q(\mathbb{R}^3)\cap D^{1,r}(\mathbb{R}^3)$, we have
			\begin{align}
			\label{2.1}	&||f||_{L^p}^p\le C||f||_{L^2}^\frac{6-p}{2}||\nabla f||_{L^2}^{\frac{3p-6}{2}},\\
			\label{2.2}	&||g||_{C(\bar{{\mathbb{R}^3}})}\le C||g||_{L^q}^{\frac{q(r-3)}{3r+q(r-3)}}||\nabla g||_{L^r}^{\frac{q(r-3)}{3r+q(r-3)}}.
			\end{align}
		\end{lemma}
		Then, by the following equalities
		\begin{align}
        &\label{2.3}\Delta\mathbf{f}=\nabla\mathrm{div}\mathbf{f}-\mathrm{curl}(\mathrm{curl}\mathbf{f}),\quad
        \mathrm{div}(\mathrm{curl}\mathbf{f})=0,\quad \mathrm{curl}(\nabla\mathbf{f})=0,
		\end{align}
		we can rewrite $(\ref{1.1})_2$ and $(\ref{1.1})_3$ as follows
		\begin{align}\label{2.4}
			&\Delta F_1 = \mathrm{div}(\rho\dot{\mathbf{u}}),\quad
			(\mu+\zeta)\Delta V_1=\mathrm{curl}\big(\rho\dot{\mathbf{u}}+2\zeta V_2\big),\\
		\label{2.5}
		&\Delta F_2 -\frac{4\zeta}{2\mu^{'}+\lambda^{'}}F_2= \mathrm{div}(\rho\dot{\mathbf{w}}),\quad
		\mu^{'}\Delta V_2-4\zeta V_2=\mathrm{curl}\big(\rho\dot{\mathbf{w}}+2\zeta V_1\big),
		\end{align}where		
		\begin{align*}
		&	\dot{\mathbf{f}}\triangleq\mathbf{f}_t+\mathbf{u}\cdot\nabla \mathbf{f},
		\quad V_1\triangleq\mathrm{curl}\mathbf{u},
		\quad F_1\triangleq(2\mu+\lambda)\mathrm{div}\mathbf{u}-P(\rho)+P(1),\\
		& V_2\triangleq\mathrm{curl}\mathbf{w},
		\quad	F_2\triangleq(2\mu^{'}+\lambda^{'})\mathrm{div}\mathbf{w}.
		\end{align*}
		Now, we state some elementary estimates which follow from the standard $L^{p}$-estimates and (\ref{2.1}). 	
		\begin{lemma}Let $(\rho,\mathbf{u},\mathbf{w})$ be the smooth solution to $(\ref{1.1})$-$(\ref{1.3})$ on $\mathbb{R}^3\times(0,T]$, for $p\in[2,6]$ it holds that
			\begin{align}
		\label{2.6}		& ||\nabla F_1||_{L^p}\le C||\rho\dot{\mathbf{u}}||_{L^p}\\
		\label{2.7}			& ||\nabla V_1||_{L^p}\le C(||\rho\dot{\mathbf{u}}||_{L^p}+||V_2||_{L^p})\\
		\label{2.8}			& ||F_1||_{L^p}\le C(||\nabla\mathbf{u}||_{L^2}+||P-P(\tilde{\rho})||_{L^2})^{\frac{6-p}{2p}}||\rho\dot{\mathbf{u}}||_{L^2}^{\frac{3p-6}{2p}}\\
		\label{2.9}			& ||V_1||_{L^p}\le C||\nabla\mathbf{u}||_{L^2}^{\frac{6-p}{2p}}(||\rho\dot{\mathbf{u}}||_{L^2}+||V_2||_{L^2})^{\frac{3p-6}{2p}}.\\
		\label{2.10}			&||\nabla\mathbf{u}||_{L^p}\le C(||F_1||_{L^p}+||V_1||_{L^p}+||P-P(\tilde{\rho})||_{L^p})	\\
		\label{2.11}		&||\nabla\mathbf{u}||_{L^p}\le||\nabla\mathbf{u}||_{L^2}^{\frac{6-p}{2p}}(||\rho\dot{\mathbf{u}}||_{L^2}+||V_2||_{L^2}+||P-P(\tilde{\rho})||_{L^6})^{\frac{3p-6}{2p}},	\\
		\label{2.12}		&||\nabla F_2||_{L^2}+||F_2||_{L^p}\le C||\rho\dot{\mathbf{w}}||_{L^2},\\
		\label{2.13}		&||\nabla F_2||_{L^p}\le C(||\rho\dot{\mathbf{w}}||_{L^p}+||\mathbf{w}||_{L^p}),\\
        \label{2.14}	    &||\nabla V_2||_{L^2}+||V_2||_{L^p}\le C(||\rho\dot{\mathbf{w}}||_{L^2}+||V_1||_{L^2}),\\
        \label{2.15}	    &||\nabla V_2||_{L^p}\le C(||\rho\dot{\mathbf{w}}||_{L^p}+||\mathbf{w}||_{L^p}+||V_1||_{L^p}),\\
		\label{2.16}		&||\nabla\mathbf{w}||_{L^p}\le||\nabla\mathbf{w}||_{L^2}^{\frac{6-p}{2p}}(||\rho\dot{\mathbf{w}}||_{L^2}+||V_1||_{L^2})^{\frac{3p-6}{2p}},          
		\end{align}		
		\end{lemma}
	{\bf Proof:} The proof of this lemma is similar to the proof of Lemma 2.3 in \cite{18}, we omit it for brevity.$\hfill\Box$	\\
	
	Next, we will use the following Zlotnik inequality to get the uniform upper bound of the density $\rho$.
		\begin{lemma}(\cite{27})
			Let the function $y$ satisfy
			\begin{align*}
				y^{'}(t)= g(y)+b^{'}(t)\quad \text{on} [0,T], \ y(0)=y^{0},
			\end{align*}
			with $g\in C(\mathbb{R})$ and $y,b\in W^{1,1}(0,T)$. If $g(\infty)=-\infty$ and
			\begin{align}
				b(t_2)-b(t_1)\le N_0+N_1(t_2-t_1)
			\end{align}
			for all $0\le t_1< t_2\le T$ with some $N_0\ge 0$ and $N_1\ge 0$, then
			\begin{align*}
				y(t)\le \max\{y^0,\bar{\xi}\}+N_0<\infty \quad \text{on}\quad [0,T],
			\end{align*} 
			where $\bar{\xi}$ is a constant such that
			\begin{align}
				g(\xi)\le-N_1\quad\text{for}\quad \xi\ge\bar{\xi}.
			\end{align}
		\end{lemma}
		Finally, the following Bcale-Kato-Majda-type inequality in \cite{28,29} will be used later to estimate $||\nabla\mathbf{u}||_{L^\infty}$ and $||\nabla\rho||_{L^2\cap L^6}$.
		\begin{lemma}\label{l2.3}
			For $3<q<\infty$, there is a constant $C(q)$ such that the following estimate holds for all $\nabla\mathbf{u}\in L^2(\mathbb{R}^3)\cap D^{1,q}(\mathbb{R}^3)$;
			\begin{align}
				||\nabla\mathbf{u}||_{L^\infty}\le C\big(||\mathrm{div}\mathbf{u}||_{L^\infty}+||\mathrm{curl}\mathbf{u}||_{L^\infty}\big)\mathrm{log}\big(e+||\nabla^2\mathbf{u}||_{L^q}\big)+C||\nabla\mathbf{u}||_{L^2}+C.
			\end{align} 
		\end{lemma}
		
		\section{A Priori Estimates}
		Let $T>0$ be a positive time and $(\rho,\mathbf{u},\mathbf{w})$ be the smooth solution to $(\ref{1.1})$-$(\ref{1.3})$ on $\mathbb{R}^3\times(0,T]$, we set 
			  \begin{align*}
			&	A_1(T)\triangleq\sup\limits_{0\le t\le T}\big(\sigma(||\nabla\mathbf{u}||_{L^2}^2+||\nabla\mathbf{w}||_{L^2}^2)\big)
				+\int_{0}^{T}\int\sigma(\rho|\mathbf{\dot{u}}|^2+\rho|\mathbf{\dot{w}}|^2)dxdt ,\\
			&	A_2(T)\triangleq\sup\limits_{0\le t\le T}\bigg(\sigma^3\int(\rho|\mathbf{\dot{w}}|^2+\rho|\mathbf{\dot{u}}|^2)dx\bigg)+\int_{0}^{T}\sigma^3\big(||\nabla\mathbf{\dot{u}}||_{L^2}^2+||\nabla\mathbf{\dot{w}}||_{L^2}^2\big)dt,\\
			&	A_3(T)\triangleq\sup\limits_{0\le t\le T}\int\rho|\mathbf{u}|^3dx,
				\end{align*}
				where	$\sigma(t)\triangleq\min\{1,t\}$.
				
			In the following, we denote the generic constant by $C>0$ depending on some known constants $\mu,\ \lambda,\ \zeta,\ \mu^{'},\ \lambda^{'},\ \gamma\ \text{and}\ M$, but independent of time $t$. Particularly, we write $C(\alpha)$ to emphasize that $C$ may depend on $\alpha$, and let $E_0 \leq 1$ without loss of generality.
				
			we give the following key a priori estimates on the solution $(\rho,\mathbf{u},\mathbf{w})$. 
       \begin{proposition}\label{p3.1} Under the conditions of Theorem $\ref{t1.1}$, if $(\rho,\mathbf{u},\mathbf{w})$ is the smooth solution to $(\ref{1.1})$-$(\ref{1.3})$ on $\mathbb{R}^3\times(0,T]$ satisfying
       	\begin{align}\label{3.1}
       		\rho\le2\bar{\rho},\quad A_1(T)+A_2(T)\le 2 E_0^{\frac{1}{2}},\quad A_3(\sigma(T))\le2E_0^{\frac{1}{4}} ,
       	\end{align}
       	then the following estimates holds
       \begin{align}\label{3.2}
       	    \rho\le\frac{7}{4}\bar{\rho},\quad A_1(T)+A_2(T)\le E_0^{\frac{1}{2}},\quad A_3(\sigma(T))\le E_0^{\frac{1}{4}} 
       \end{align}
       provieded $E_0\le \epsilon$, where $\epsilon$ is a positive constant depending on $\mu,\ \lambda,\ \zeta,\ \mu^{'},\ \lambda^{'},\ \gamma\ \text{and}\ M$.
       \end{proposition}
       {\bf Proof:} Propositon \ref{p3.1} is proved by Lemma \ref{l3.4}-\ref{l3.6}.$\hfill\Box$	\\

	 We give the standard energy estimate for $(\rho,\mathbf{u},\mathbf{w})$ at first		
	\begin{lemma}\label{l3.1} Let $(\rho,\mathbf{u},\mathbf{w})$ be the smooth solution to $(\ref{1.1})$-$(\ref{1.3})$ on $\mathbb{R}^3\times(0,T]$ with $0\le\rho(x,t)\le 2\bar{\rho}$, it holds that
	\begin{align}\label{3.3}
	\begin{aligned}
	&\sup_{0\let\le T} \int\big(\frac{1}{2}\rho|{\mathbf{u}}|^2+\frac{1}{2}\rho|{\mathbf{w}}|^2+\frac{1}{2}|\mathrm{\mathbf{B}}|^2+G(\rho)\big)dx+\int\big(\mu|\nabla\mathrm{\mathbf{u}}|^2+(\lambda+\mu)(\mathrm{div}\mathrm{\mathbf{u}})^2\\&+\int_{0}^{T}\mu^{'}|\nabla\mathrm{\mathbf{w}}|^2+(\lambda^{'}+\mu^{'})(\mathrm{div}\mathrm{\mathbf{w}})^2+\zeta|2\mathbf{w}-\mathrm{curl}\mathbf{u}|^2\big)dxdt \le E_0.	
	\end{aligned}
     \end{align}
     \end{lemma}
     {\bf Proof:} Multiplying $(\ref{1.1})_1$-$(\ref{1.1})_3$ by $G^{'}(\rho)$, $\mathbf{u}$ and $\mathbf{w}$ respectively, then integrating by parts, one get the $(\ref{3.3})$ directly.$\hfill\Box$
     
    \begin{lemma}\label{l3.2}
    	 Let $(\rho,\mathbf{u},\mathbf{w})$ be the smooth solution to $(\ref{1.1})$-$(\ref{1.3})$ on $\mathbb{R}^3\times(0,T]$ with $0\le\rho(x,t)\le 2\bar{\rho}$, it holds that	   	
	\begin{align}\label{3.4}
	 A_1(T)\le CE_0+C\int_{0}^{T}\sigma(||\nabla\mathbf{u}||_{L^3}^3+||\nabla\mathbf{w}||_{L^3}^3)dt
	 \end{align}
	 and
	 \begin{align}\label{3.5}
	   \begin{aligned}
	   & A_2(T)+\int_{0}^{T}\sigma^3||\mathrm{curl}\mathbf{\dot{u}}-2\mathbf{\dot{w}}||_{L^2}^2dt\\\le& CE_0+A_1(T)+C\int_{0}^{T}\sigma^3||\nabla\mathbf{u}||_{L^2}||\nabla\mathbf{w}||_{L^3}^3dt+C\int_{0}^{T}\sigma^{3}\big(||\nabla\mathbf{w}||_{L^4}^4+||\nabla\mathbf{u}||_{L^4}^4\big)dt.  
	   \end{aligned}
	 \end{align}		   	
	\end{lemma}	
	{\bf Proof:} For $m\ge0$, multiplying the $(\ref{1.1})_2$ and $(\ref{1.1})_3$ by $\sigma^m\dot{\mathbf{u}}$ and $\sigma^m\dot{\mathbf{u}}$ respectively, integrating the resulting equalities over $\mathbb{R}^3$ and summing them up leads to 
	\begin{align}\label{3.6}
	\begin{aligned}
		\sigma^m\int\big(\rho|\mathbf{\dot{u}}|^2+\rho|\mathbf{\dot{w}}|^2\big)dx=&-\sigma^m\int\nabla P(\rho)\cdot\mathbf{\dot{u}}dx+(\mu+\zeta)\sigma^m\int\Delta\mathbf{u}\cdot\mathbf{\dot{u}}dx\\&+(\mu+\lambda-\zeta)\sigma^m\int\nabla\mathrm{div}\mathbf{u}\cdot\mathbf{\dot{u}}dx+2\zeta \int\mathrm{curl}{\mathbf{w}}\cdot\mathbf{\dot{u}}dx\\&-4\zeta\sigma^m\int\mathbf{w}\cdot\mathbf{\dot{w}}dx+\mu^{'}\sigma^m\int\Delta\mathbf{w}\cdot\mathbf{\dot{w}}dx\\&+(\mu^{'}+\lambda^{'})\sigma^m\int\nabla\mathrm{div}\mathbf{w}\cdot\mathbf{\dot{w}}dx+2\zeta \int\mathrm{curl}{\mathbf{u}}\cdot\mathbf{\dot{w}}dx\\
		:=& \sum_{j=1}^{8}I_j.
		\end{aligned}
    \end{align}
    Similar to the proof of Lemma 3.2 in \cite{18}, we give the estimates of $I_1$, $I_2+I_3$ and $I_6+I_7$ as follows,
	\begin{align*}
	I_1\le\bigg(\int\sigma^m\mathrm{div}\mathbf{u}(P(\rho)-P(\tilde{\rho}))dx\bigg)_t+C(\bar{\rho})||\nabla\mathbf{u}||_{L^2}^2+C(\bar{\rho})m^2\sigma^{2(m-1)}\sigma^{'}E_0,
	\end{align*}	
	\begin{align*}
    \begin{aligned}
	I_2+I_3
	 \le& -\bigg(\frac{\mu+\zeta}{2}\sigma^m||\nabla\mathbf{u}||_{L^2}^2+\frac{\mu+\lambda-\zeta}{2}\sigma^m||\mathrm{div}\mathbf{u}||_{L^2}^2\bigg)_t\\
	&+Cm\sigma^{m-1}\sigma^{'}||\nabla\mathbf{u}||_{L^2}^2+C\sigma^m||\nabla\mathbf{u}||_{L^3}^3\\
	\le&-\big(\frac{\mu}{2}\sigma^m||\nabla\mathbf{u}||_{L^2}^2+\frac{\mu+\lambda}{2}\sigma^m||\mathrm{div}\mathbf{u}||_{L^2}^2+\frac{\zeta}{2}\sigma^m||\mathrm{curl}\mathbf{u}||_{L^2}^2\big)_t\\&
	+Cm\sigma^{m-1}\sigma^{'}||\nabla\mathbf{u}||_{L^2}^2+C\sigma^m||\nabla\mathbf{u}||_{L^3}^3,
	\end{aligned}
	\end{align*}
and
	\begin{align*}
	\begin{aligned}
	 I_6+I_7
	\le &-\bigg(\frac{\mu^{'}}{2}\sigma^m||\nabla\mathbf{w}||_{L^2}^2+\frac{\mu^{'}+\lambda^{'}}{2}\sigma^m||\mathrm{div}\mathbf{w}||_{L^2}^2\bigg)_t\\
	&+Cm\sigma^{m-1}\sigma^{'}||\nabla\mathbf{w}||_{L^2}^2+C\sigma^m(||\nabla\mathbf{u}||_{L^3}^3+||\nabla\mathbf{w}||_{L^3}^3).\\
	\end{aligned}
	\end{align*}
By H\"older inequality and $(\ref{2.1})$, then integrating by parts leads to
   \begin{align*}
   \begin{aligned}
    I_4+I_8=& 4\zeta\sigma^m\int\mathrm{curl}\mathbf{u}_t\cdot\mathbf{w}dx+2\zeta\sigma^m\int\mathbf{u}\cdot\nabla\mathbf{u}\cdot\mathrm{curl}\mathbf{w}dx+2\zeta\sigma^m\int\mathbf{u}\cdot\nabla\mathbf{w}\cdot\mathrm{curl}\mathbf{u}dx\\
   \le& 4\zeta\sigma^m\int\mathrm{curl}\mathbf{u}_t\cdot\mathbf{w}dx+C\sigma^m||\mathbf{u}||_{L^6}||\nabla\mathbf{u}||_{L^3}||\nabla\mathbf{w}||_{L^2}\\
   \le& 4\zeta\sigma^m\int\mathrm{curl}\mathbf{u}_t\cdot\mathbf{w}dx+ C\sigma^m||\nabla\mathbf{u}||_{L^3}^3+C\sigma^m(||\nabla\mathbf{u}||_{L^2}^3+||\nabla\mathbf{w}||_{L^2}^3)
   \end{aligned}
   \end{align*}
   and
   \begin{align*}
   \begin{aligned}
   I_5=&-2\zeta\bigg(\sigma^m\int|\mathbf{w}|^2dx\bigg)_t+2\zeta m\sigma^{m-1}\sigma^{'}\int|\mathbf{w}|^2dx-4\zeta\sigma^m\int\mathbf{w}\cdot\mathbf{u}\cdot\nabla\mathbf{w}dx\\
   \le& -2\zeta\bigg(\sigma^m\int|\mathbf{w}|^2dx\bigg)_t+2\zeta m\sigma^{m-1}\sigma^{'}||\mathbf{w}||_{L^2}^2+C\sigma^m(||\nabla\mathbf{u}||_{L^2}^2+||\nabla\mathbf{w}||_{L^2}^3)
   \end{aligned}
   \end{align*}
  Inserting the estimates of $I_j$ $(j=1,\cdots,8)$ into 
  $(\ref{3.6})$, one has
   \begin{align}\label{3.7}
   \begin{aligned}
   &\bigg(\sigma^m\big(\frac{\mu}{2}||\nabla\mathbf{u}||_{L^2}^2+\frac{\mu+\lambda}{2}||\mathrm{div}\mathbf{u}||_{L^2}^2+\frac{\mu^{'}}{2}||\nabla\mathbf{w}||_{L^2}^2+\frac{\mu^{'}+\lambda^{'}}{2}||\mathrm{div}\mathbf{w}||_{L^2}^2\big)\bigg)_t\\&+\frac{\zeta}{2}\big(\sigma^m||2\mathbf{w}-\mathrm{curl}\mathbf{u}||_{L^2}^2\big)_t+\sigma^m\int(\rho|\mathbf{\dot{u}}|^2+\rho|\mathbf{\dot{w}}|^2)dx  \\
   \le& C m\sigma^{m-1}\sigma^{'}\big(\frac{\mu}{2}||\nabla\mathbf{u}||_{L^2}^2+\frac{\mu+\lambda}{2}||\mathrm{div}\mathbf{u}||_{L^2}^2+\frac{\mu^{'}}{2}||\nabla\mathbf{w}||_{L^2}^2+\frac{\mu^{'}+\lambda^{'}}{2}||\mathrm{div}\mathbf{w}||_{L^2}^2\big)\\&+Cm\sigma^{m-1}\sigma^{'}\frac{\zeta}{2}||2\mathbf{w}-\mathrm{curl}\mathbf{u}||_{L^2}^2+C\sigma^m(||\nabla\mathbf{u}||_{L^3}^3+||\nabla\mathbf{w}||_{L^3}^3). 
   \end{aligned}
   \end{align}
   Let $m=1$ in $(\ref{3.7})$, and integrating resulting inequality over $(0,T)$, it follows from $(\ref{3.3})$ that $(\ref{3.4})$.
 
  Next, for $m\ge0$, operating $\sigma^m\dot{\mathbf{u}}^j[\partial_t+\mathrm{div}(\mathbf{u}\cdot)]$ and $\sigma^m\dot{\mathbf{w}}^j[\partial_t+\mathrm{div}(\mathbf{u}\cdot)]$ on $(\ref{1.1})_2$ and $(\ref{1.1})_3$ respectively, summing the resulting equalities up, we get
   \begin{align}\label{3.8}
   \begin{aligned}
   &\bigg(\frac{\sigma^m}{2}\int\big(\rho|\mathbf{\dot{u}}|^2+\rho|\mathbf{\dot{w}}|^2\big)dx\bigg)_t-\frac{m}{2}\sigma^{m-1}\sigma^{'}\int\big(\rho|\mathbf{\dot{u}}|^2+\rho|\mathbf{\dot{w}}|^2\big)dx\\
   =&-\int\sigma^m\dot{\mathbf{u}}^j\big(\partial_jP_t+\mathrm{div}(\partial_jP\mathbf{u})\big)dx\\&+\int\sigma^m(\mu+\zeta)\dot{\mathbf{u}}^j\big(\Delta\mathbf{u}_t^j+\mathrm{div}(\mathbf{u}\Delta\mathbf{u}^j)\big)dx\\
   &+(\mu+\lambda-\zeta)\int\sigma^m\dot{\mathbf{u}}^j\big(\partial_j\mathrm{div}\mathbf{u}_t+\mathrm{div}(\mathbf{u}\partial_j\mathrm{div}\mathbf{u})\big)dx\\&+2\zeta\int\sigma^m\dot{\mathbf{u}}^j\bigg((\mathrm{curl}\mathbf{w})_t^j+\mathrm{div}\big(\mathbf{u}(\mathrm{curl}\mathbf{w})^j\big)\bigg)dx\\
   &-4\zeta\int\sigma^m\dot{\mathbf{w}}^j\big(\mathbf{w}_t+\mathrm{div}(\mathbf{u}\mathbf{w}^j)\big)dx\\&+\mu^{'}\int\sigma^m\dot{\mathbf{w}}^j\big(\Delta\mathbf{w}_t^j+\mathrm{div}(\mathbf{u}\Delta\mathbf{w}^j)\big)dx\\
   &+(\mu^{'}+\zeta^{'})\int\sigma^m\dot{\mathbf{w}}^j\big(\partial_j\mathrm{div}\mathbf{w}_t+\mathrm{div}(\mathbf{u}\partial_j\mathrm{div}\mathbf{w})\big)dx\\&+2\zeta\int\sigma^m\dot{\mathbf{w}}^j\bigg((\mathrm{curl}\mathbf{u})_t^j+\mathrm{div}\big(\mathbf{u}(\mathrm{curl}\mathbf{u})^j\big)\bigg)dx\\
   =& \sum_{i=1}^{8}J_i.
   \end{aligned}
   \end{align}    
   Similar to the proof of Lemma 3.2 in \cite{18}, we give the estimates of $J_1+J_2+J_3$ and $J_6+J_7$ as follows,
   \begin{align*}
    \begin{aligned}
      J_1+J_2+J_3\le&-\big(\frac{\mu}{2}\sigma^m||\nabla\mathbf{\dot{u}}||_{L^2}^2+\frac{\mu+\lambda}{2}\sigma^m||\mathrm{div}\mathbf{\dot{u}}||_{L^2}^2+{\zeta}\sigma^m||\mathrm{curl}\mathbf{\dot{u}}||_{L^2}^2\big)\\&
      +C(\bar{\rho})\sigma^{m}||\nabla\mathbf{u}||_{L^2}^2+C\sigma^m||\nabla\mathbf{u}||_{L^4}^4 
      \end{aligned}
   \end{align*}
   and
   \begin{align*}
   \begin{aligned}
      J_6+J_7\le&-\big(\frac{\mu^{'}}{2}\sigma^m||\nabla\mathbf{\dot{w}}||_{L^2}^2+\frac{\mu^{'}+\lambda^{'}}{2}\sigma^m||\mathrm{div}\mathbf{\dot{w}}||_{L^2}^2\big)\\&
      +C\sigma^{m}(||\nabla\mathbf{u}||_{L^4}^4+||\nabla\mathbf{w}||_{L^4}^4).   
    \end{aligned}
   \end{align*}
   It follows from H\"older inequality, (\ref{3.1}) and (\ref{3.3}) that 
   \begin{align*}
   \begin{aligned}
   	J_4&=2\zeta\int\sigma^m\mathrm{curl}\dot{\mathbf{u}}\cdot\mathbf{w}_tdx-2\zeta\int\sigma^m\nabla\dot{\mathbf{u}}^j\cdot\big(\mathbf{u}\mathrm{curl}\mathbf{w}^j\big)dx\\
   	&\le2\zeta\int\sigma^m\mathrm{curl}\dot{\mathbf{u}}\cdot\mathbf{\dot{w}}dx+\sigma^m||\nabla\dot{\mathbf{u}}||_{L^2}||\mathbf{u}||_{L^6}||\mathrm{curl}\mathbf{w}||_{L^3}\\
   	&\le2\zeta\int\sigma^m\mathrm{curl}\dot{\mathbf{u}}\cdot\mathbf{\dot{w}}dx+\frac{\mu}{8}\sigma^m||\nabla\dot{\mathbf{u}}||_{L^2}^2+C\sigma^m||\nabla\mathbf{u}||_{L^2}^4+C\sigma^m||\nabla\mathbf{u}||_{L^2}||\nabla\mathbf{w}||_{L^3}^3\\
   	&\le2\zeta\int\sigma^m\mathrm{curl}\dot{\mathbf{u}}\cdot\mathbf{\dot{w}}dx+\frac{\mu}{8}\sigma^m||\nabla\dot{\mathbf{u}}||_{L^2}^2+CE_0^\frac{1}{2}\sigma^{m-1}||\nabla\mathbf{u}||_{L^2}^2+C\sigma^m||\nabla\mathbf{u}||_{L^2}||\nabla\mathbf{w}||_{L^3}^3,
   	\end{aligned}
   \end{align*}
      \begin{align*}
      	\begin{aligned}
      	J_5\le 4\zeta||\mathbf{\dot{w}}||_{L^2}^2+C\sigma^m||\nabla{\mathbf{u}}||_{L^2}||\mathbf{w}||_{L^6}||\mathbf{w}||_{L^3}\le 4\zeta||\mathbf{\dot{w}}||_{L^2}^2+ C\sigma^{m-1}E_0^{\frac{1}{2}}(||{\mathbf{w}}||_{L^2}^2+||\nabla{\mathbf{w}}||_{L^2}^2)
      	\end{aligned}
      \end{align*}
and
   \begin{align*}
   \begin{aligned}
   	 J_8 &=2\zeta\int\sigma^m\dot{\mathbf{w}}\cdot\mathrm{curl}\mathbf{u}_tdx+2\zeta\int\sigma^m\dot{\mathbf{w}}^j\cdot\big(\mathbf{u}\nabla\mathrm{curl}\mathbf{w}^j\big)dx+2\zeta\int\sigma^m\dot{\mathbf{w}}^j(\mathrm{div}\mathbf{u}\cdot\mathrm{curl}\mathbf{w}^j)dx\\ 	 &\le2\zeta\int\sigma^m\mathrm{curl}\dot{\mathbf{u}}\cdot\mathbf{\dot{w}}dx+C\sigma^m||\dot{\mathbf{w}}||_{L^6}||\nabla\mathbf{u}||_{L^2}||\mathrm{curl}\mathbf{w}||_{L^3}\\
   	 &\le2\zeta\int\sigma^m\mathrm{curl}\dot{\mathbf{u}}\cdot\mathbf{\dot{w}}dx+\frac{\mu^{'}}{8}\sigma^m||\nabla\dot{\mathbf{w}}||_{L^2}^2+C\sigma^m||\nabla\mathbf{u}||_{L^2}^4+C\sigma^m||\nabla\mathbf{u}||_{L^2}||\nabla\mathbf{w}||_{L^3}^3\\
   	 &\le2\zeta\int\sigma^m\mathrm{curl}\dot{\mathbf{u}}\cdot\mathbf{\dot{w}}dx+\frac{\mu^{'}}{8}\sigma^m||\nabla\dot{\mathbf{w}}||_{L^2}^2+CE_0^\frac{1}{2}\sigma^{m-1}||\nabla\mathbf{u}||_{L^2}^2+C\sigma^m||\nabla\mathbf{u}||_{L^2}||\nabla\mathbf{w}||_{L^3}^3.
   	 \end{aligned}
   \end{align*}
   Inserting the estimates of $J_i$ $(i=1,\cdots,8)$ into 
   $(\ref{3.8})$ and integrating resulting inequality over $(0,T)$, it follows from $(\ref{3.3})$ that
   \begin{align}\label{3.9}
   \begin{aligned}
   &\frac{\sigma^m}{2}\int\big(\rho|\mathbf{\dot{w}}|^2+\rho|\mathbf{\dot{u}}|^2\big)dx
   +\int_{0}^{T}\sigma^m\big(\frac{\mu}{2}||\nabla\mathbf{\dot{u}}||_{L^2}^2+\frac{\mu+\lambda}{2}||\mathrm{div}\mathbf{\dot{u}}||_{L^2}^2\\&+\frac{\mu^{'}}{2}||\nabla\mathbf{\dot{w}}||_{L^2}^2+\frac{\mu^{'}+\lambda^{'}}{2}||\mathrm{div}\mathbf{\dot{w}}||_{L^2}^2+\zeta||\mathrm{curl}\mathbf{\dot{u}}-2\mathbf{\dot{w}}||_{L^2}^2\big)dt\\
   \le&CE_0+\frac{m}{2}\int_{0}^{T}\sigma^{m-1}\sigma^{'}\int\big(\rho|\mathbf{\dot{u}}|^2+\rho|\mathbf{\dot{w}}|^2\big)dx\\
   &+C\int_{0}^{T}\sigma^m||\nabla\mathbf{u}||_{L^2}||\nabla\mathbf{w}||_{L^3}^3dt+C\int_{0}^{T}\sigma^{m}\big(||\nabla\mathbf{w}||_{L^4}^4+||\nabla\mathbf{u}||_{L^4}^4\big)dt.
   \end{aligned}
   \end{align}
      Let $m=3$ in $(\ref{3.9})$, one gets $(\ref{3.5})$.Thus, we finish the proof of Lemma \ref{l3.2}.$\hfill\Box$
   
 \begin{lemma}\label{l3.3} 
 	Let $(\rho,\mathbf{u},\mathbf{w})$ be the smooth solution to $(\ref{1.1})$-$(\ref{1.3})$ on $\mathbb{R}^3\times(0,T]$ with $0\le\rho(x,t)\le 2\bar{\rho}$, it holds that	
		\begin{align}\label{3.10}
	  &  \sup\limits_{0\le t\le \sigma({T})}\big(||\nabla\mathbf{u}||_{L^2}^2+||\nabla\mathbf{w}||_{L^2}^2\big)+\int_{0}^{\sigma({T})}\int\rho(|\mathbf{\dot{u}}|^2+|\mathbf{\dot{w}}|^2)dxdt 
		\le C(M)
	\end{align}and
	\begin{align}
\label{3.11}
	  &  \sup\limits_{0\le t\le \sigma({T})} t\int\rho\big(|\mathbf{\dot{w}}|^2+|\mathbf{\dot{u}}|^2\big)dx
	    +\int_{0}^{\sigma(T)}t \big(||\nabla\mathbf{\dot{u}}||_{L^2}^2+
	    ||\nabla\mathbf{\dot{w}}||_{L^2}^2\big)dt
	    \le C(M) 	.
	    \end{align}	
 \end{lemma}  
 {\bf Proof:} For $0\le t\le \sigma(T)$, multiplying $(\ref{1.2})_2$ and $(\ref{1.2})_3$ by $\mathbf{u}_t$ and $\mathbf{w}_t$, integrating the resulting equalities by parts over $\mathbb{R}^3$, we obtain by H\"older inequality, $(\ref{2.1})$, $(\ref{2.6})$, $(\ref{2.11})$ and $(\ref{2.16})$ that
       \begin{align}\label{3.12}
       \begin{aligned}
       &\frac{1}{2}\frac{d}{dt}\bigg({\mu}||\nabla\mathbf{u}||_{L^2}^2+({\mu+\lambda})||\mathrm{div}\mathbf{u}||_{L^2}^2+{\zeta}||\mathrm{curl}\mathbf{u}||_{L^2}^2\bigg)+ \int\rho|\mathbf{\dot{\mathbf{u}}}|^2dx-2\zeta\int\mathrm{curl}{\mathbf{u}}_t\cdot\mathbf{w}dx \\
       =&\frac{d}{dt}\int\big(P(\rho)-P(1)\big)\mathrm{div}\mathbf{u}dx+\int\rho\dot{\mathbf{u}}(\mathbf{u}\cdot\nabla\mathbf{u})dx\\&+\int\mathrm{div}\mathbf{u}\mathrm{div}\big(\big(P(\rho)-P(1)\big)\mathbf{u}\big)dx+\int\big((\gamma-1)P(\rho)+P(1)\big)\mathrm{div}\mathbf{u}\mathrm{div}\mathbf{u}dx\\
       \le& \frac{d}{dt}\int\big(P(\rho)-P(1)\big)\mathrm{div}\mathbf{u}dx+C\bigg(\int\rho|\dot{\mathbf{u}}|^2dx\bigg)^\frac{1}{2}\bigg(\int\rho|\mathbf{u}|^3\bigg)^{\frac{1}{3}}||\nabla\mathbf{u}||_{L^6}\\&-\frac{1}{2\mu+\lambda}\int\nabla\big(F_1+(P(\rho)-P(1))\big)\cdot\big(P(\rho)-P(1)\big)\mathbf{u}dx\\
       \le& \frac{d}{dt}\int\big(P(\rho)-P(1)\big)\mathrm{div}\mathbf{u}dx\\&+CE_0^{\frac{1}{12}}||\rho^\frac{1}{2}\dot{\mathbf{u}}||_{L^2}\big(||\rho^\frac{1}{2}\dot{\mathbf{u}}||_{L^2}+||\mathrm{curl}\mathbf{w}||_{L^2}+||P(\rho)-P(1)||_{L^6}\big)\\&
       +C||\rho^\frac{1}{2}\dot{\mathbf{u}}||_{L^2}||P(\rho)-P(1)||_{L^3}||{\mathbf{u}}||_{L^6}+C||P(\rho)-P(1)||_{L^4}^2||\nabla{\mathbf{u}}||_{L^2}\\
       \le& \frac{d}{dt}\int\big(P(\rho)-P(1)\big)\mathrm{div}\mathbf{u}dx+C E_0^{\frac{1}{12}}\int\rho|\dot{\mathbf{u}}|^2dx+C\big(||\nabla\mathbf{u}||_{L^2}^2+||\mathrm{curl}\mathbf{w}||_{L^2}^2\big)\\&+C||P(\rho)-P(1)||_{L^4}^4+C||P(\rho)-P(1)||_{L^6}^2
       \end{aligned}
       \end{align}
       and
      \begin{align}\label{3.13}
      \begin{aligned}
      &\frac{1}{2}\frac{d}{dt}\big(\mu^{'}||\nabla{\mathbf{w}}||_{L^2}^2+(\mu^{'}+\lambda^{'})||\mathrm{div}\mathbf{w}||_{L^2}^2+4\zeta||\mathbf{w}||_{L^2}^2\big)+\int\rho |\dot{\mathbf{w}}|^2dx-2\zeta\int\mathrm{curl}\mathbf{u}{\mathbf{w}}_tdx\\
      =&\int\rho\dot{\mathbf{w}}\cdot\mathbf{u}\nabla\mathbf{w}dx\\
      \le& C\bigg(\int\rho|\dot{\mathbf{w}}|^2dx\bigg)^{\frac{1}{2}}\bigg(\int\rho\mathbf{u}^3dx\bigg)^{\frac{1}{3}}||\nabla\mathbf{w}||_{L^6}\\
      \le& C\bigg(\int\rho|\dot{\mathbf{w}}|^2dx\bigg)^{\frac{1}{2}}\bigg(\int\rho\mathbf{u}^3dx\bigg)^{\frac{1}{3}}(||\rho\dot{\mathbf{w}}||_{L^2}+||\mathrm{curl}\mathbf{u}||_{L^2})\\
      \le& CE_0^{\frac{1}{12}}\bigg(\int\rho|\dot{\mathbf{w}}|^2dx+||\mathrm{curl}\mathbf{u}||_{L^2}^2\bigg)
      \end{aligned}
      \end{align}
      Summing $(\ref{3.12})$ and $(\ref{3.13})$ up, integrating over $(0,\sigma(T))$, gives $(\ref{3.10})$.
       
       Let $T=\sigma(T)$ and $m=1$ in $(\ref{3.9})$, Using H\"older inequality, (\ref{2.11}), (\ref{2.16}), (\ref{3.1}), (\ref{3.3}) and (\ref{3.10}), one has
          \begin{align}
          \begin{aligned}
          &\sup\limits_{0\le t\le\sigma(T)}t\int\big(\rho|\mathbf{\dot{w}}|^2+\rho|\mathbf{\dot{u}}|^2\big)dx
          +\int_{0}^{\sigma(T)}t\big({\mu}||\nabla\mathbf{\dot{u}}||_{L^2}^2+({\mu+\lambda})||\mathrm{div}\mathbf{\dot{u}}||_{L^2}^2\\&+{\mu^{'}}||\nabla\mathbf{\dot{w}}||_{L^2}^2+({\mu^{'}+\lambda^{'}})||\mathrm{div}\mathbf{\dot{w}}||_{L^2}^2+2\zeta||2\mathbf{\dot{w}}-\mathrm{curl}\mathbf{\dot{u}}||_{L^2}^2\big)dt\\
          \le& C\int_{0}^{\sigma(T)}\int\big(\rho|\mathbf{\dot{u}}|^2+\rho|\mathbf{\dot{w}}|^2\big)dx+C\int_{0}^{\sigma(T)}\big(||\nabla{\mathbf{u}}||_{L^2}^2+||\nabla{\mathbf{w}}||_{L^2}^2\big)dt\\
          &+C\int_{0}^{\sigma(T)}\sigma||\nabla\mathbf{u}||_{L^2}||\nabla\mathbf{w}||_{L^3}^3dt+C\int_{0}^{\sigma(T)}\sigma\big(||\nabla\mathbf{w}||_{L^4}^4+||\nabla\mathbf{u}||_{L^4}^4\big)dt.\\
         \le & C(M)+C(M)\int_{0}^{\sigma(T)}\big(\sigma||\nabla\mathbf{w}||_{L^6}^2+\sigma(||\nabla\mathbf{w}||_{L^6}^3+||\nabla\mathbf{u}||_{L^6}^3)\big)dt\\
         \le& C (M)+C(M)\int_{0}^{\sigma(T)}\sigma\big(||\rho\dot{\mathbf{w}}||_{L^2}^2+||\mathrm{curl}\mathbf{w}||_{L^2}^2+||P-P(\tilde{\rho})||_{L^6}^2\big)dt\\&+\int_{0}^{\sigma(T)}\sigma\big(||\rho\dot{\mathbf{w}}||_{L^2}^3+||\rho\dot{\mathbf{u}}||_{L^2}^3+||\mathrm{curl}\mathbf{w}||_{L^2}^3+||\mathrm{curl}\mathbf{u}||_{L^2}^3+||P-P(\tilde{\rho})||_{L^6}^3\big)dt      \\   
         \le& C(M).
          \end{aligned}
          \end{align}   
          Thus, we finish the proof of Lemma \ref{l3.3}.$\hfill\Box$
  \begin{lemma}\label{l3.4} Let $(\rho,\mathbf{u},\mathbf{w})$ be the smooth solution to $(\ref{1.1})$-$(\ref{1.3})$ on $\mathbb{R}^3\times(0,T]$ with $0\le\rho(x,t)\le 2\bar{\rho}$, it holds that	
  	    \begin{align}
  	    	\sup\limits_{0\le t\le T}\int\rho |\mathbf{u}|^3dx\le E_0^{\frac{1}{4}},
  	    \end{align}
  	    provided that $E_0\le\big(\frac{1}{C(M)}\big)^2$.
  \end{lemma}        
  {\bf Proof:} Multiplying $(\ref{1.1})_2$ by $3|\mathbf{u}|\mathbf{u}$ and integrating the resulting equation over $\mathbb{R}^3$, it follows from $(\ref{2.11})$, $(\ref{3.1})$, $(\ref{3.3})$ and $(\ref{3.10})$ that
   \begin{align*}
  \begin{aligned}
  	&\int\rho|\mathbf{u}|^3dx\\
  	\le&\int\rho_0|\mathbf{u}_0|^3dx+ C\int_{0}^{\sigma(T)}\int|\mathbf{u}||\nabla\mathbf{u}|^2dxdt+C\int_{0}^{\sigma(T)}\int|P(\rho)-P(1)||\mathbf{u}||\nabla\mathbf{u}|dxdt\\&+C\int_{0}^{\sigma(T)}\int|\mathbf{w}||\mathbf{u}||\nabla\mathbf{u}|dxdt\\
  	\le&\int\rho_0|\mathbf{u}_0|^3dx+ C \int_{0}^{\sigma(T)}||\mathbf{u}||_{L^6}||\nabla\mathbf{u}||_{L^2}^\frac{3}{2}||\nabla\mathbf{u}||_{L^6}^{\frac{1}{2}}dt\\&+C\int_{0}^{\sigma(T)}||P(\rho)-P(1)||_{L^3}||\mathbf{u}||_{L^6}||\nabla\mathbf{u}||_{L^2}dt+C\int_{0}^{\sigma(T)}||\mathbf{w}||_{L^3}||\mathbf{u}||_{L^6}||\nabla\mathbf{u}||_{L^2}dt\\
  	\le& C\bigg(\int\rho_0|\mathbf{u}_0|^2dx\bigg)^{\frac{3}{4}}\bigg(\int|\mathbf{u}_0|^6dx\bigg)^\frac{1}{4}+ CE_0^{\frac{1}{3}}\int_{0}^{\sigma(T)}||\nabla\mathbf{u}||_{L^2}^{2}dt\\&C\int_{0}^{\sigma(T)}||\nabla\mathbf{u}||_{L^2}^{\frac{5}{2}}\big(||\rho\dot{\mathbf{u}}||_{L^2}+||\mathrm{curl}\mathbf{w}||_{L^2}+||P-P(\tilde{\rho})||_{L^6}\big)^{\frac{1}{2}}dt\\
  	\le& C(M)E_0^{\frac{3}{4}}+CE_0+C(M)\bigg(\int_{0}^{\sigma({T})}||\nabla\mathbf{u}||_{L^2}^2dt\bigg)^\frac{1}{4}\bigg(\int_{0}^{\sigma({T})}||\nabla\mathbf{u}||_{L^2}^4dt\bigg)^\frac{1}{2}\bigg(\int_{0}^{\sigma({T})}||\rho\dot{\mathbf{u}}||_{L^2}^2dt\bigg)^\frac{1}{4}\\
  	\le& C(M)E_0^{\frac{3}{4}}.
  \end{aligned}  
    \end{align*}  	   Let $E_0\le\big(\frac{1}{C(M)}\big)^2$, we complete the proof of this lemma. $\hfill\Box$
    \begin{lemma}\label{l3.5}Let $(\rho,\mathbf{u},\mathbf{w})$ be the smooth solution to $(\ref{1.1})$-$(\ref{1.3})$ on $\mathbb{R}^3\times(0,T]$ with $0\le\rho(x,t)\le 2\bar{\rho}$, it holds that
    	\begin{align*}
    	A_1(T)+A_2(T)\le E_0^{\frac{1}{2}},
    	\end{align*} provided that $E_0\le\min\big\{\big(\frac{1}{2C(M)}\big)^8,\big(\frac{1}{2C}\big)^2\big\}$.
    \end{lemma}
      	    provided that $E_0\le\big(\frac{1}{C(M)}\big)^8$.
    {\bf Proof:} Adding (\ref{3.4}) and (\ref{3.5}) up, we obtain
    	\begin{align}\label{3.14}\begin{aligned}
    	& A_1(T)+A_2(T)\\\le& CE_0+C\int_{0}^{T}\sigma\big(||\nabla\mathbf{u}||_{L^3}^3+||\nabla\mathbf{w}||_{L^3}^3\big)dt+C\int_{0}^{T}\sigma^3\big(||\nabla\mathbf{u}||_{L^4}^4+||\nabla\mathbf{w}||_{L^4}^4\big)dt.   	\end{aligned}
    	\end{align}
    	It follows from $(\ref{2.16})$, $(\ref{3.1})$ and $(\ref{3.3})$ that
    	\begin{align}\label{3.15}\begin{aligned}
        \int_{0}^{T}\sigma^3||\nabla\mathbf{w}||_{L^4}^4dt
    	\le	&\int_{0}^{T}\sigma^3||\nabla\mathbf{w}||_{L^2}\big(||\rho\dot{\mathbf{w}}||_{L^2}+||\nabla\mathbf{u}||_{L^2}\big)^3dt\\
    	\le& 
    	\sup\limits_{0\le t\le T}\sigma^3\big(||\nabla\mathbf{w}||_{L^2}^2+||\nabla\mathbf{u}||_{L^2}^2\big)
    	\int_{0}^{T}||\nabla\mathbf{u}||_{L^2}^2dt\\
    	&+\sup\limits_{0\le t \le T}(\sigma^\frac{1}{2}||\nabla\mathbf{w}||_{L^2})(\sigma^\frac{3}{2} ||\sqrt{\rho}\dot{\mathbf{u}}||_{L^2})\int_{0}^{T}\sigma\int\rho|\dot{\mathbf{u}}|^2dxdt\\
    \le	& CE_0.	
    	\end{aligned}
    	\end{align}
    	Due to $(\ref{2.10})$,
        \begin{align}\label{3.16}
        \begin{aligned}
        &\int_{0}^{T}\sigma^3||\nabla\mathbf{u}||_{L^4}^4dt
		        \le C\int_{0}^{T}\sigma^3(||F_1||_{L^4}^4+||V_1||_{L^4}^4)dt+C\int_{0}^{T}\sigma^3||P-P(\tilde{\rho})||_{L^4}^4dt\\
        \end{aligned}
        \end{align}
        By $(\ref{2.8})$, $(\ref{2.9})$, $(\ref{3.1})$ and $(\ref{3.3})$, we get
        \begin{align}\label{3.17}
        \begin{aligned}
        &\int_{0}^{T}\sigma^3(||F_1||_{L^4}^4+||V_1||_{L^4}^4)dt\\
        \le& C\int_{0}^{T}\sigma^3(||\nabla\mathbf{u}||_{L^2}+||P-P(\tilde{\rho})||_{L^2})||\rho\dot{\mathbf{u}}||_{L^2}^3dt+C\int_{0}^{T}\sigma^3||\nabla\mathbf{u}||_{L^2}||V_2||_{L^2}^3dt\\
        \le& C\sup\limits_{0 \le t \le T}\big(\sigma^\frac{3}{2}||\rho\dot{\mathbf{u}}||_{L^2}(||\nabla\mathbf{u}||_{L^2}+E_0^\frac{1}{2})\big)\int_{0}^{T}\sigma||\rho^\frac{1}{2}\dot{\mathbf{u}}||_{L^2}^2dt\\&+C\sup\limits_{0 \le t \le T}\sigma^3\big(||\nabla\mathbf{u}||_{L^2}^2+||\nabla\mathbf{w}||_{L^2}^2\big)\int_{0}^{T}||\nabla\mathbf{w}||_{L^2}^2dt   \\
        \le& CE_0 .     
        \end{aligned}
        \end{align}
        Similar to the proof of Lemma $\ref{3.5}$ in $\cite{18}$, one has
        \begin{align}\label{3.18}
        \int_{0}^{T}\sigma^3||P(\rho)-P(1)||_{L^4}^4dt\le CE_0.
        \end{align}
        Applying $(\ref{3.16})$-$(\ref{3.18})$, we obtain
         \begin{align}\label{3.19}
         \int_{0}^{T}\sigma^3\big(||\nabla\mathbf{u}||_{L^4}^4+||P(\rho)-P(1)||_{L^4}^4\big)dt\le CE_0.
         \end{align}
         By $(\ref{3.15})$ and $(\ref{3.19})$, we obtain
    	\begin{align}\label{3.20}
    	\begin{aligned}
    	&\int_{\sigma(T)}^{T}\sigma\big(||\nabla\mathbf{w}||_{L^3}^3+||\nabla\mathbf{u}||_{L^3}^3\big)dt\\
    	\le&\int_{\sigma(T)}^{T}\big(||\nabla\mathbf{w}||_{L^2}^2+||\nabla\mathbf{w}||_{L^4}^4+||\nabla\mathbf{u}||_{L^2}^2+||\nabla\mathbf{u}||_{L^4}^4\big)dt\\
    	\le& CE_0.    
    	\end{aligned}
    	\end{align}
    	With the help of (\ref{2.8}), (\ref{2.9}), (\ref{2.16}), (\ref{3.1}), (\ref{3.10}) and (\ref{3.11}), we get
        \begin{align}\label{3.21}
        \begin{aligned}
        &\int_{0}^{\sigma(T)}\sigma\big(||\nabla\mathbf{u}||_{L^3}^3+||\nabla\mathbf{w}||_{L^3}^3\big)dt\\
        \le& C\int_{0}^{\sigma(T)}\sigma||\nabla\mathbf{u}||_{L^2}^\frac{3}{2}\big(||\rho\dot{\mathbf{u}}||_{L^2}^\frac{3}{2}+||\nabla\mathbf{w}||_{L^2}^\frac{3}{2}+||P(\rho)-P(1)||_{L^2}^\frac{3}{2}\big)dt\\&+C\int_{0}^{\sigma(T)}\sigma||\nabla\mathbf{w}||_{L^2}^\frac{3}{2}||\rho\dot{\mathbf{w}}||_{L^2}^\frac{3}{2}dt\\
        \le&\sup\limits_{0 \le t \le T}(\sigma||\nabla\mathbf{u}||_{L^2}^2)^\frac{1}{4}\int_{0}^{\sigma(T)}||\nabla\mathbf{u}||_{L^2}
        (\sigma||\rho\dot{\mathbf{u}}||_{L^2}^2)^\frac{3}{4}dt\\&+\sup\limits_{0 \le t \le T}(\sigma||\nabla\mathbf{w}||_{L^2}^2)^\frac{1}{4}\int_{0}^{\sigma(T)}||\nabla\mathbf{w}||_{L^2}
        (\sigma||\rho\dot{\mathbf{w}}||_{L^2}^2)^\frac{3}{4}dt+C E_0\\
        \le& C(M)A_1(T)^\frac{1}{4}\big(E_0+A_1(T)\big)+C E_0\\
        \le&  C(M)E_0^\frac{5}{8}.   
        \end{aligned}
        \end{align} 
	 Let $E_0\le\min\big\{\big(\frac{1}{2C(M)}\big)^8,\big(\frac{1}{2C}\big)^2\big\}$, we finish this lemma by inserting $(\ref{3.17})$ and $(\ref{3.19})$-$(\ref{3.21})$ into $(\ref{3.16})$. $\hfill\Box$
    \begin{lemma}\label{l3.6}Let $(\rho,\mathbf{u},\mathbf{w})$ be the smooth solution to $(\ref{1.1})$-$(\ref{1.3})$ on $\mathbb{R}^3\times(0,T]$ with $0\le\rho(x,t)\le 2\bar{\rho}$, it holds that
    	\begin{align}\label{3.22}
    		\sup\limits_{0\le t\le T}||\rho(t)||_{L^\infty}\le\frac{7\bar{\rho}}{4},   		
    	\end{align}
    	provided that $E_0\le\min\big\{\big(\frac{\bar{\rho}}{2C(M)}\big)^8,\big(\frac{\bar{\rho}}{4C(M)}\big)^\frac{3}{2}\big\}$
    \end{lemma}    
        {\bf Proof:} By $(\ref{1.1})_1$, we have
        \begin{align*}
        	D_t\rho=g(\rho)+b^{'}(t)
        \end{align*}
        where
        \begin{align*}
        	D_t\rho=\rho_t+\mathbf{u}\cdot\nabla\rho,\quad g(\rho)=-\frac{a\rho}{2\mu+\lambda}(\rho^\gamma-\tilde{\rho}^\gamma),\quad b(t)=\frac{1}{2\mu+\lambda}\int_{0}^{t}\rho F_1dt.
        \end{align*}
        For $t\in[0,\sigma(T)]$, it follows from for all $0\le t_1\le t_2\le\sigma(T)$ that
        \begin{align*}
        	|b(t_2)-b(t_1)|&\le C\int_{0}^{\sigma(T)}||\rho F_1||_{L^{\infty}}dt\\
        	&\le\int_{0}^{T}||F_1||_{L^6}^\frac{1}{2}||\nabla F_1||_{L^6}^\frac{1}{2}dt\\
        	&\le\int_{0}^{T}||\rho\dot{\mathbf{u}}||_{L^2}^\frac{1}{2}||\nabla\dot{\mathbf{u}}||_{L^2}^\frac{1}{2}dt\\
        	&\le C(M)\int_{0}^{T}t^{-\frac{1}{4}}||\rho\dot{\mathbf{u}}||_{L^2}^\frac{1}{2}(t||\nabla\dot{\mathbf{u}}||_{L^2}^2)^\frac{1}{4}dt\\
        	&\le C(M)\bigg(\int_{0}^{T}t^{-\frac{1}{3}}||\rho\dot{\mathbf{u}}||_{L^2}^\frac{2}{3}dt\bigg)^{\frac{3}{4}}\\
        	&\le C(M)(A_1(\sigma(T)))^\frac{1}{4}\\
        	&\le C(M)E_0^\frac{1}{8}.
        \end{align*}By Lemma\ref{l2.3}, let  
        \begin{align*}
        	N_1=0,\quad N_0=C(M)E_0^\frac{1}{8},\quad \bar{\xi}=\tilde{\rho}.
        \end{align*}Then, one has
        \begin{align*}
        	g(\xi)=-\frac{a\xi}{2\mu+\lambda}(\xi^\gamma-\tilde{\rho}^\gamma)\le-N_1=0\quad \text{for all}\ {\xi\ge\bar{\xi}=\tilde{\rho}}.
        \end{align*}Thus, Let $E_0\le \big(\frac{\bar{\rho}}{2C(M)}\big)^8$, one has
                \begin{align}\label{p1}
        	\sup\limits_{0\let\le\sigma(T)}||\rho||_{L^\infty}\le\max\big\{\bar{\rho},\tilde{\rho}\big\}+N_0\le\bar{\rho}+C(M)E_0^\frac{1}{8}\le\frac{3\bar{\rho}}{2}.
        \end{align}
         For $t\in[\sigma(T),T]$, it follows from for all $\sigma(T)\le t_1\le t_2\le T$ that
         \begin{align*}
         |b(t_2)-b(t_1)|&\le C\int_{t_1}^{t_2}||\rho F_1||_{L^{\infty}}dt\\
         &\le \frac{a}{2\mu+\lambda}(t_2-t_1)+C\int_{\sigma(T)}^{T}||F_1||_{L^{\infty}}^\frac{8}{3}dt\\
                 	&\le \frac{a}{2\mu+\lambda}(t_2-t_1)+C\int_{0}^{T}||F_1||_{L^2}^\frac{2}{3}||\nabla F_1||_{L^6}^2dt\\
         &\le\frac{a}{2\mu+\lambda}(t_2-t_1)+CE_0^\frac{1}{6}\int_{\sigma(T)}^{T}||\nabla\dot{\mathbf{u}}||_{L^2}^2dt\\
         &\le\frac{a}{2\mu+\lambda}(t_2-t_1)+ C(M)E_0^\frac{2}{3}.
         \end{align*}By Lemma\ref{l2.3}, let  
                 \begin{align*}
                 N_1=\frac{a}{2\mu+\lambda},\quad N_0=C(M)E_0^\frac{2}{3},\quad \bar{\xi}=\tilde{\rho}+1.
                 \end{align*}Then, one has
                 \begin{align*}
                 g(\xi)=-\frac{a\xi}{2\mu+\lambda}(\xi^{\gamma}-\tilde{\rho}^{\gamma})\le-N_1=\frac{a}{2\mu+\lambda}\quad \text{for all} \ \xi\ge\tilde{\rho}+1.
                 \end{align*}Thus, Let $E_0\le \big(\frac{\bar{\rho}}{4C(M)}\big)^\frac{3}{2}$, we obtain
                 \begin{align}\label{p2}
                 \sup\limits_{0\le t\le\sigma(T)}||\rho||_{L^\infty}\le\max\big\{\frac{3}{2}\bar{\rho},\tilde{\rho}+1\big\}+N_0\le\frac{3\bar{\rho}}{2}+CE_0^
                 	\frac{2}{3}\le\frac{7\bar{\rho}}{4}.
                 \end{align}
                 $\hfill\Box$
                 
                 From now on, we will always assume that the positive constant $C$ may depend on $T$, $||\rho_0g_i||_{L^2}\ (i=1,2)$, $||\nabla \mathbf{u}_0||_{H^1}$, $||\nabla \mathbf{w}_0||_{H^1}$, $||\rho_0-\tilde{\rho}||_{H^2}$, $||P(\rho_0)-P(\tilde{\rho})||_{H^2}$ besides $\mu,\ \lambda,\ \zeta,\ \mu^{'},\ \lambda^{'},\ \gamma\ \text{and}\ M$. 
                 
                 Next, we will give the higher order estimates for smooth solution of $(\ref{1.1})$-$(\ref{1.3})$. 
     \begin{lemma}\label{l3.7}Under the conditions of Theorem $\ref{t1.1}$, for any given
     	$T>0$, it holds that
		\begin{align}\label{3.23}
		&\sup\limits_{0\le t\le T}\big(||\nabla\mathbf{u}||_{L^2}^2+||\nabla\mathbf{w}||_{L^2}^2\big)+\int_{0}^{T}\big(||\rho^\frac{1}{2}\mathbf{\dot{u}}||_{L^2}^2+||\rho^\frac{1}{2}\mathbf{\dot{w}}||_{L^2}^2+||\mathrm{curl}\mathbf{\dot{u}}-2\mathbf{\dot{w}}||_{L^2}^2\big)dt \le C,\\
		\label{3.24}&\sup\limits_{0\le t\le T}\big(||\rho^\frac{1}{2}\mathbf{\dot{u}}||_{L^2}^2+||\rho^\frac{1}{2}\mathbf{\dot{w}}||_{L^2}^2\big)+\int_{0}^{T}\big(||\nabla\mathbf{\dot{u}}||_{L^2}^2+||\nabla\mathbf{\dot{w}}||_{L^2}^2\big)dt	\le C.
		\end{align}
     \end{lemma}
     {\bf Proof:}  Using (\ref{3.1}) and (\ref{3.10}), we get (\ref{3.23}) directly.
     
     Taking $m=0$ in (\ref{3.9}), it follows from (\ref{2.6}), (\ref{2.7}), (\ref{2.10}), (\ref{2.12}),(\ref{2.14}),(\ref{2.16}) and (\ref{3.23}) that
     \begin{align}\label{3.25}
     \begin{aligned}
     &\frac{d}{dt}\big(||\rho^\frac{1}{2}\mathbf{\dot{u}}||_{L^2}^2+||\rho^\frac{1}{2}\mathbf{\dot{w}}||_{L^2}^2\big)
     +\big(\frac{\mu}{2}||\nabla\mathbf{\dot{u}}||_{L^2}^2+\frac{\mu+\lambda}{2}||\mathrm{div}\mathbf{\dot{u}}||_{L^2}^2+\frac{\mu^{'}}{2}||\nabla\mathbf{\dot{w}}||_{L^2}^2\\&+\frac{\mu^{'}+\lambda^{'}}{2}||\mathrm{div}\mathbf{\dot{w}}||_{L^2}^2+\zeta||\mathrm{curl}\mathbf{\dot{u}}-2\mathbf{\dot{w}}||_{L^2}^2\big)\\
     \le& C
     +C||\nabla\mathbf{u}||_{L^2}||\nabla\mathbf{w}||_{L^3}^3+C\big(||\nabla\mathbf{w}||_{L^4}^4+||\nabla\mathbf{u}||_{L^4}^4\big)\\
     \le& C\big(||\nabla\mathbf{u}||_{L^6}^3+||\nabla\mathbf{w}||_{L^6}^3\big)+C\\
     \le& C\big(||F_1||_{L^6}^3+||V_1||_{L^6}^3+||P(\rho)-P(1)||_{L^6}^3+||F_2||_{L^6}^3+||V_2||_{L^6}^3\big)+C\\
     \le& C\big(||\nabla F_1||_{L^2}^3+||\nabla V_1||_{L^2}^3+||\nabla F_2||_{L^2}^3+||\nabla V_2||_{L^6}^3\big)+C\\
     \le& C\big(||\rho\mathbf{\dot{u}}||_{L^2}^3+||\rho\mathbf{\dot{w}}||_{L^2}^3\big)+C\\
     \le& C\big(||\rho^\frac{1}{2}\mathbf{\dot{u}}||_{L^2}^4+||\rho^\frac{1}{2}\mathbf{\dot{w}}||_{L^2}^4\big)+C.
     \end{aligned}
     \end{align}  
      By (\ref{3.25}) and Granwall inequality, one obtains (\ref{3.24}). Thus, we complete the proof of this lemma. $\hfill\Box$
    
    \begin{lemma}\label{l3.8}Under the conditions of Theorem $\ref{t1.1}$, for any given
    	$T>0$, it holds that
     \begin{align}\label{3.26}
     \sup\limits_{0\le t\le T}\big(||\nabla\rho||_{L^2}+||\nabla\rho||_{L^6}+||\nabla\mathbf{u}||_{H^1}\big)+\int_{0}^{T}||\nabla\mathbf{u}||_{L^\infty}dt\le C
     \end{align}
    \end{lemma}
    {\bf Proof:}  
    Based on Lemma \ref{l3.7} and the Beale-Kato-Majda type inequality (see Lemma \ref{l2.3}),
    we can establish the following estimates similar to \cite{18}.$\hfill\Box$.
   \begin{lemma}\label{l3.9}Under the conditions of Theorem $\ref{t1.1}$, for any given
   	$T>0$, it holds that
   \begin{align}
  \label{3.27} &\sup\limits_{0\le t\le T}\int\big(\rho|\mathbf{u}_t|^2+\rho|\mathbf{w}_t|^2\big)dx+\int_{0}^{T}\big(||\nabla\mathbf{u}_t||_{L^2}^2+||\nabla\mathbf{w}_t||_{L^2}^2\big)dt\le C,\\
   \label{3.28}&   \sup\limits_{0\le t\le T}\int\big(||\rho||_{H^2}+
      ||P||_{H^2}\big)dx+\int_{0}^{T}\big(||\nabla^2\mathbf{u}||_{H^1}+||\nabla^2\mathbf{w}||_{H^1}\big)dt\le C,\\
  \label{3.29}& \sup\limits_{0\le t\le T}\int\big(||\rho_t||_{H^1}+
   ||P_t||_{H^1}\big)dx+\int_{0}^{T}\big(||\rho_{tt}||_{L^2}+
   ||P_{tt}||_{L^2}\big)dt\le C.
   \end{align}   
   \end{lemma}
   {\bf Proof:} See the proofs of (3.62), (3.63), and (3.66) in \cite{18}.
   \begin{lemma}\label{l3.10}
   	Under the conditions of Theorem $\ref{t1.1}$, for any given
   	$T>0$, it holds that
   \begin{align}\label{3.30}
   \begin{aligned}
   &\sup\limits_{0\le t\le T}\sigma\big(||\nabla \mathrm{\mathbf{u}_t}||_{L^2}^2+||\nabla \mathrm{\mathbf{w}_t}||_{L^2}^2+||\nabla\mathbf{u}||_{H^2}^2+||\nabla\mathbf{w}||_{H^2}^2\big)\\&+\int_{0}^{T}\sigma\big(||\rho^{\frac{1}{2}} \mathrm{\mathbf{u}_{tt}}||_{L^2}^2+||\rho^{\frac{1}{2}} \mathrm{\mathbf{w}_{tt}}||_{L^2}^2+||\nabla \mathrm{\mathbf{u}_{t}}||_{H^1}^2+||\nabla \mathrm{\mathbf{w}_{t}}||_{H^1}^2\big)dt\le C
   \end{aligned}	
   \end{align}
   \end{lemma}
   {\bf Proof:} Differentiating $(\ref{1.1})_2$ and  $(\ref{1.1})_3$ with respect to $t$ gives
   \begin{align}
   \label{3.31}	\rho \mathbf{u}_{tt}-(\mu+\zeta)\Delta \mathrm{\mathbf{u}_t}-(\mu+\lambda-\zeta)\nabla \mathrm{div}\mathbf{u}_t-2\zeta\mathrm{curl}\mathbf{w}_t&=-\rho_t \mathbf{u}_t- (\rho\mathbf{u}\nabla\mathbf{u})_t-\nabla P_t,\\
   	\label{3.32} \rho \mathbf{w}_{tt}-\mu^{'}\Delta \mathrm{\mathbf{w}_t}-(\mu^{'}+\lambda^{'})\nabla \mathrm{div}\mathbf{w}_t+4\zeta\mathbf{w}_t-2\zeta\mathrm{curl}\mathbf{u}_t&=-\rho_t \mathbf{w}_t- (\rho\mathbf{w}\nabla\mathbf{w})_t.
   \end{align}
   Multiplying $(\ref{3.31})$ and $(\ref{3.32})$ by $\mathbf{u}_{tt}$ and $\mathbf{w}_{tt}$, respectively, and integrating the resulting equation by parts, it follows from Lemma {\ref{2.1} and \ref{l3.7}-\ref{l3.9} and (\ref{1.1}) that
   \begin{align}\label{3.33}
   \begin{aligned}
   	&\frac{1}{2}\frac{d}{dt}\big(\mu||\nabla \mathrm{\mathbf{u}_t}||_{L^2}^2+(\mu+\lambda) ||\mathrm{div}\mathrm{\mathbf{u}_t}||_{L^2}^2+\zeta ||\mathrm{curl}\mathrm{\mathbf{u}_t}||_{L^2}^2\big)+||\rho^\frac{1}{2}\mathbf{u}_{tt}||_{L^2}^2-2\zeta\int\mathrm{curl}\mathbf{u}_{t}\cdot\mathbf{w}_{tt}\\
   	=&-\frac{d}{dt}\int\big[\frac{1}{2}\rho_t|\mathbf{u}_t|^2+(\rho_t\mathbf{u}\cdot\nabla\mathbf{u}+\nabla P_t)\cdot\mathbf{u}_t\big]dx+\frac{1}{2}\int\rho_{tt}|\mathbf{u}_t|^2dx-\int P_{tt}\mathrm{div}\mathbf{u}_tdx\\
   	&+\int\big(\rho_{tt}\mathbf{u}\nabla\mathbf{u}+\rho_t\mathbf{u}_t\cdot\nabla\mathbf{u}+\rho_t\mathbf{u}\cdot\nabla\mathbf{u}_t\big)\cdot\mathbf{u}_tdx-\int(\rho\mathbf{u}_t\cdot\nabla\mathbf{u}+\rho\mathbf{u}\cdot\nabla\mathbf{u}_t)\cdot\mathbf{u}_{tt}dx\\
   	\le&-\int\big(\rho\mathbf{u}\cdot\nabla\mathbf{u}_t\cdot\mathbf{u}_t+\rho_t\mathbf{u}\cdot\nabla\mathbf{u}\cdot\mathbf{u}_t-P_t\mathrm{div}\mathbf{u}_t\big)dx+\frac{1}{2}||\rho^\frac{1}{2}\mathbf{u}_{tt}||_{L^2}^2\\
   	&+C\big(1+||\rho_{tt}||_{L^2}^2+||P_{tt}||_{L^2}^2+||\nabla\mathbf{u}_{t}||_{L^2}^2\big) \\
   \le&\frac{\mu}{4}||\nabla\mathbf{u}_t||_{L^2}^2+\frac{1}{2}||\rho^\frac{1}{2}\mathbf{u}_{tt}||_{L^2}^2+C\big(1+||\rho_{tt}||_{L^2}^2+||P_{tt}||_{L^2}^2+||\nabla\mathbf{u}_{t}||_{L^2}^2\big)
   	\end{aligned}
   \end{align}
  and
   \begin{align}\label{3.34}
   \begin{aligned}
   &\frac{1}{2}\frac{d}{dt}\big(\mu^{'}||\nabla \mathrm{\mathbf{w}_t}||_{L^2}^2+(\mu^{'}+\lambda^{'}) ||\mathrm{div}\mathrm{\mathbf{w}_t}||_{L^2}^2+4\zeta ||\mathrm{\mathbf{w}_t}||_{L^2}^2\big)+||\rho^\frac{1}{2}\mathbf{w}_{tt}||_{L^2}^2-2\zeta\int\mathrm{curl}\mathbf{u}_{tt}\cdot\mathbf{w}_t\\
   =&-\frac{d}{dt}\int\big(\frac{1}{2}\rho_t|\mathbf{w}_t|^2+\rho_t\mathbf{u}\cdot\nabla\mathbf{w}\cdot\mathbf{w}_t\big)dx+\frac{1}{2}\int\rho_{tt}|\mathbf{w}_t|^2dx\\
   &+\int\big(\rho_{tt}\mathbf{u}\nabla\mathbf{w}+\rho_t\mathbf{u}_t\cdot\nabla\mathbf{w}+\rho_t\mathbf{u}\cdot\nabla\mathbf{w}_t\big)\cdot\mathbf{w}_tdx-\int(\rho\mathbf{u}_t\cdot\nabla\mathbf{w}+\rho\mathbf{u}\cdot\nabla\mathbf{w}_t)\cdot\mathbf{u}_{tt}dx\\
   \le&-\int\big(\rho\mathbf{u}\cdot\nabla\mathbf{w}_t\cdot\mathbf{w}_t+\rho_t\mathbf{u}\cdot\nabla\mathbf{w}\cdot\mathbf{w}_t\big)dx \\
   & +\frac{1}{2}||\rho^\frac{1}{2}\mathbf{w}_{tt}||_{L^2}^2 +C\big(1+||\rho_{tt}||_{L^2}^2+||P_{tt}||_{L^2}^2+||\nabla\mathbf{u}_{t}||_{L^2}^2+||\nabla\mathbf{w}_{t}||_{L^2}^2\big)\\
   \le&\frac{\mu^{'}}{4}||\nabla\mathbf{w}_t||_{L^2}^2+\frac{1}{2}||\rho^\frac{1}{2}\mathbf{w}_{tt}||_{L^2}^2+C\big(1+||\rho_{tt}||_{L^2}^2+||P_{tt}||_{L^2}^2+||\nabla\mathbf{u}_{t}||_{L^2}^2+||\nabla\mathbf{w}_{t}||_{L^2}^2\big)
   \end{aligned}
   \end{align}
Adding (\ref{3.33}) and (\ref{3.34}) up and multiplying resulting inequality by $\sigma(t)$, integrating it over $[0,T]$, by(\ref{2.1}), (\ref{3.25}) and (\ref{3.27})-(\ref{3.30}), we obtain
   \begin{align}\label{3.35}
   \begin{aligned}
   &\sup\limits_{0\le t\le T}\sigma\big(||\nabla \mathrm{\mathbf{u}_t}||_{L^2}^2+||\nabla \mathrm{\mathbf{w}_t}||_{L^2}^2\big)+\int_{0}^{T}\sigma\big(||\rho^{\frac{1}{2}} \mathrm{\mathbf{u}_{tt}}||_{L^2}^2+||\rho^{\frac{1}{2}} \mathrm{\mathbf{w}_{tt}}||_{L^2}^2\big)dt\\\le& C\int_{0}^{T}\sigma^{'}\int|\mathrm{curl}\mathbf{u}_t-2\mathbf{w}_t|^2dxdt+C\\
   \le& C\int_{0}^{T}||\mathrm{curl}\mathbf{\dot{u}}-2\mathbf{\dot{w}}||_{L^2}^2dt+C\int_{0}^{T}\int\big|\mathrm{curl}(\mathbf{u}\nabla\mathbf{u})-2\mathbf{u}\nabla\mathbf{w}\big|^2dxdt+C\\
   \le& C.
   \end{aligned}
   \end{align}
   By (\ref{3.27})-(\ref{3.30}), we infer from $(1.1)_2$-$(1.1)_3$ and the standard $L^2$-estimate of the elliptic system that
   \begin{align}\label{3.36}
   \begin{aligned}
   	&||\nabla\mathbf{u}||_{H^2}^2+||\nabla\mathbf{w}||_{H^2}^2\\
   	\le& C+C\big(||\nabla(\rho\mathbf{u}_t+\rho\mathbf{u}\cdot\nabla\mathbf{u})||_{L^2}^2+||\nabla^2P||_{L^2}+||\nabla\mathbf{w}||_{H^1}\big)\\&+C\big(||\nabla(\rho\mathbf{w}_t+\rho\mathbf{w}\cdot\nabla\mathbf{w})||_{L^2}^2+||\nabla^2\mathbf{u}||_{L^2}\big)\\
   	\le& C\big(1+||\nabla\mathbf{u}_t||_{L^2}+||\nabla\mathbf{w}_t||_{L^2}\big), \end{aligned}
   \end{align}     
 
   Similarly,
   \begin{align}\label{3.37}
   	\begin{aligned}
   	&||\nabla^2\mathbf{u}_t||_{L^2}+||\nabla^2\mathbf{w}_t||_{L^2}\\
   	\le&||\big(\rho\mathbf{u}_t+\rho\mathbf{u}\nabla\mathbf{u}+\nabla P(\rho)+\mathrm{curl}\mathbf{w}+\rho\mathbf{w}_{t}+\rho\mathbf{u}\nabla\mathbf{w}+\mathbf{w}+\mathrm{curl}\mathbf{u}\big)_t||_{L^2}\\
   	\le& C\big(1+||\rho^\frac{1}{2}\mathbf{u}_{tt}||_{L^2}+||\rho^\frac{1}{2}\mathbf{w}_{tt}||_{L^2}+||\nabla\mathbf{u}_t||_{L^2}+||\nabla\mathbf{w}_t||_{L^2} \big)
   	\end{aligned}
   \end{align}
      (\ref{3.36}) and (\ref{3.37}) combined with (\ref{3.35}), gives
      \begin{align}\label{3.38}
      \sup\limits_{0\le t \le T}\sigma\big(||\nabla\mathbf{u}||_{H^2}^2+||\nabla\mathbf{w}||_{H^2}^2\big)+\int_{0}^{T}\sigma\big(||\nabla^2\mathbf{u}_t||_{L^2}^2+||\nabla^2\mathbf{w}_t||_{L^2}^2\big)dt
      \le C.
      \end{align}  
   By (\ref{3.35}) and (\ref{3.38}), we finish the proof of this lemma.$\hfill\Box$
   \begin{lemma}\label{l3.11}
   For $3<q<6$, under the conditions of Theorem $\ref{t1.1}$, for any given
   $T>0$, it holds that
   \begin{align}\label{3.39}
   \begin{aligned}
   &\sup\limits_{0\le t \le T}\big(||\nabla\rho||_{w^{1,q}}+||\nabla P||_{w^{1,q}}\big)\\&+\int_{0}^{T}\big(||\nabla\mathbf{u}_t||_{L^q}^{p_0}+||\nabla^2\mathbf{u}||_{W^{1,q}}^{p_0}+||\nabla\mathbf{w}_t||_{L^q}^{p_0}+||\nabla^2\mathbf{w}||_{W^{1,q}}^{p_0}\big)dt\le C,  
   \end{aligned}
   \end{align}
   where\begin{align*}
   1\le p_0\le \frac{4q}{5q-6}\in (1,2).
   \end{align*}
   \end{lemma}
  {\bf Proof:}Based on Lemma \ref{l3.7}-\ref{l3.10}, the proof of this lemma is similar to Lemma 4.5 in \cite{19}.$\hfill\Box$
    \begin{lemma}\label{l3.12}
        	Under the conditions of Theorem $\ref{t1.1}$, for any given
        	$T>0$, it holds that
        	\begin{align}\label{3.40}\begin{aligned}
        	&\sup\limits_{0\le t\le T}\sigma\big(||\rho^{\frac{1}{2}} \mathrm{\mathbf{u}_{tt}}||_{L^2}+||\rho^{\frac{1}{2}} \mathrm{\mathbf{w}_{tt}}||_{L^2}+||\nabla^2\mathbf{u}||_{W^{1,q}}+||\nabla\mathbf{w}||_{W^{1,q}}+||\nabla^2 \mathrm{\mathbf{u}_{t}}||_{L^2}+||\nabla^2 \mathrm{\mathbf{w}_{t}}||_{L^2}\big)\\&+\int_{0}^{T}\sigma^2\big(||\nabla \mathrm{\mathbf{u}_{tt}}||_{L^2}^2+||\nabla \mathrm{\mathbf{w}_{tt}}||_{L^2}^2\big)dt\le C .       	
        	\end{aligned}	
        	\end{align}
    \end{lemma}
     {\bf Proof:}Based on Lemma \ref{l3.7}-\ref{l3.11}, the proof of this lemma is similar to Lemma 4.6 in \cite{19}.$\hfill\Box$
 \section{Proof of Theorem{\ref{t1.1}}}
 We first need the following local existence theorem of classical solution of (\ref{1.1})--(\ref{1.3}) with
 large initial data.
 \begin{lemma}\label{l4.1}
 	Assume that the initial data $(\rho_0,\mathbf{u}_0,\mathbf{w}_0)$ satisfy the conditions
 	$(\ref{1.4})$, $(\ref{1.5})$ of Theorem \ref{t1.1}. Then there exist a positive time $T_0 >0$ and a unique
 	classical solution $(\rho_0,\mathbf{u}_0,\mathbf{w}_0)$ of (\ref{1.1})–(\ref{1.3}) on $\mathbb{R}^3\times(0,T_0]$, satisfying $\rho_0\ge0$, and for any $\tau\in(0,T_0)$,
 				\begin{align}\label{4.1}
 				\begin{cases}
 				&(\rho-\tilde{\rho},P(\rho)-P(\tilde{\rho})\in C([0,T_0];H^2\cap W^{2,q}),\\
 				&\mathbf{u}\in C([0,T_0];D^1\cap D^2)\cap L^\infty ([\tau,T_0];D^2\cap D^{3,q}),\\
 				&\mathbf{u}_t\in L^\infty([\tau,T_0];D^1\cap D^2)\cap H^1([\tau,T_0];D^1),\\
 				&\mathbf{w}\in C([0,T_0];D^1\cap D^2)\cap L^\infty ([\tau,T_0];D^2\cap D^{3,q}),\\
 				&\mathbf{w}_t\in L^\infty([\tau,T_0];D^1\cap D^2)\cap H^1([\tau,T_0];D^1).	
 				\end{cases}
 				\end{align}
 \end{lemma}
 {\bf Proof:} By Galerkin’s approximation method, we construct the approximate solutions $\mathbf{u}^m$ to the momentum equation. Then using this approximate $\mathbf{u}^m$ and the equations of conservation mass and  microrotational velocity, we get $\rho^m$ and $\mathbf{w}^m$.  Similar to that of the compressible Navier–Stokes equations (see, for example, [32]), the existence of a smooth approximate solution $(\rho^m,\mathbf{u}^m,\mathbf{w}^m)$ follows from the fixed point theorem. In order to prove the convergence for the approximate solutions and to obtain a smooth solution of $(\ref{1.1})$–$(\ref{1.3})$, we give some uniform a priori estimates for $(\rho^m,\mathbf{u}^m,\mathbf{w}^m)$.
 As in [], one can show that there exists a small $T_0>0$, independent of $m$ and the lower bound of density, so that
 \begin{align}\label{4.2}
 \begin{aligned}
 	&\sup\limits_{0\le t\le T_0}\big(||\rho^m-\tilde{\rho}||_{H^1\cap W^{1,q}}+ ||\sqrt{\rho^m}\mathbf{u_t}||_{L^2}^2+||\sqrt{\rho^m}\mathbf{w_t}||_{L^2}^2)+||\nabla\mathbf{u}^m||_{H^1}^2+||\nabla\mathbf{w}^m||_{H^1}^2\big)\\&+\int_{0}^{T_0}\big(||\mathbf{u}^m||_{D^{2,q}}^2+||\mathbf{w}^m||_{D^{2,q}}^2\big)dt+\int_{0}^{T_0}\big(||\nabla\mathbf{u}_t^m||_{L^2}^2+||\nabla\mathbf{w}_t^m||_{L^2}^2\big)dt
 	\le\tilde{C}. \end{aligned}
 \end{align}
 By (\ref{4.2}) we can establish the estimates on the higher-order derivatives of the solution $(\rho^m,\mathbf{u}^m,\mathbf{w}^m)$ in the same way as carried out in Lemmas ${\ref{l3.9}}$--${\ref{l3.12}}$ as follows,
    \begin{align}\label{4.3}
    \begin{aligned}
    &\sup\limits_{0\le t \le T_0}\big(||\nabla\rho^m||_{w^{1,q}}+||\nabla P^m||_{w^{1,q}}\big)\\&+\int_{0}^{T}\big(||\nabla\mathbf{u}^m_t||_{L^q}^{p_0}+||\nabla^2\mathbf{u}^m||_{W^{1,q}}^{p_0}+||\nabla\mathbf{w}^m_t||_{L^q}^{p_0}+||\nabla^2\mathbf{w}^m||_{W^{1,q}}^{p_0}\big)dt\le \tilde{C},  
    \end{aligned}
    \end{align}
    and
    \begin{align}\label{4.4}\begin{aligned}
     &\sup\limits_{\tau\le t\le T_0}\big(||\sqrt{\rho^m}\mathrm{\mathbf{u^m}_{tt}}||_{L^2}+||\sqrt{\rho^m} \mathrm{\mathbf{w^m}_{tt}}||_{L^2}+||\nabla^2\mathbf{u^m}||_{W^{1,q}}+||\nabla\mathbf{w^m}||_{W^{1,q}}\\&+||\nabla^2 \mathrm{\mathbf{u^m}_{t}}||_{H^1}+||\nabla \mathrm{\mathbf{w^m}_{t}}||_{H^1}\big)+\int_{\tau}^{T_0}\big(||\nabla \mathrm{\mathbf{u^m}_{tt}}||_{L^2}^2+||\nabla \mathrm{\mathbf{w^m}_{tt}}||_{L^2}^2\big)dt\le\tilde{C} .       	
      \end{aligned}	
    \end{align}
    {\bf Proof of Theorem \ref{1.1}:} Using lemma \ref{l4.1} and by the Proposition \ref{p3.1} and the estimates in Lemmas \ref{l3.9}-\ref{l3.12}, we can extend the local classical solution $(\rho,\mathbf{u},\mathbf{w})$ to all time, and one has
    \begin{align}\label{4.5}
     T^*=\mathrm{sup}\{\ T|\ (\ref{3.2})\ holds\ \}=\infty.
    \end{align} 
    Then, we can obtain that $(\rho,\mathbf{u},\mathbf{w})$ is continuous in $t$, especially, for $q\in(3,6)$, it shows that
    \begin{align}\label{4.6}
    \begin{cases}
    	&(\rho-\tilde{\rho},P(\rho)-P(\tilde{\rho}))\in C([0,T];D^2\cap D^{2,q}),\\
        &(\rho-\tilde{\rho},P(\rho)-P(\tilde{\rho}))\in C([0,T];H^2\cap W^{2,q}),\\
        &(\mathbf{u},\mathbf{w})\in C([0,T];D^1\cap D^2).
    \end{cases}
    \end{align} The proof of $(\ref{4.5})$ and $(\ref{4.6})$ is similar to \cite{19}, we omit it here.
    
    The large-time behavior $(\ref{1.10})$  can be obtained in similar arguments as used in \cite{18,19}. Thus, the proof of Theorem \ref{1.1} is finished.$\hfill\Box$
       \section*{Acknowledgement}
       

\begin{thebibliography}{10}    
       	\bibitem{16} Y. Amirat, K. Hamdache, Global weak solutions to a ferrofluid flow model, Mathematical Methods in the Applied Sciences., 31(2) (2008), 123--151.
       	\bibitem{26} O. A. Ladyzenskaja, V. A. Solonnikov, N. N. Ural'ceva, Linear and quasi-linear equations of parabolic type, American Mathematical Soc., (Vol. 23) (1968).  	
       	\bibitem{28} J. T. Beale, T. Kato, A. Majda, Remarks on the breakdown of smooth solutions for the
       	3-D Euler equations, Commun. Math. Phys., 94 (1984), 61-–66.	 
       	\bibitem{1} S.C. Cowin, Polar fluids, Phys. Fluids 11 (1968) 1919-–1927.
        \bibitem{4} M.T. Chen, Global strong solutions for the viscous, micropolar, compressible flow, J. Partial Differ. Equ., 24 (2011)
        158–-164. 
        \bibitem{11} M.T. Chen, Unique solvability of compressible micropolar viscous fluids, Boundary Value Problems, 1 (2012), 32.
       	\bibitem{25} Y. Cho, H. Kim,  Existence results for viscous polytropic fluids with vacuum, Journal of Differential Equations., 228(2) (2006), 377--411.        	
       	\bibitem{24} H. J. Choe, H. Kim, Strong solutions of the Navier–Stokes equations for isentropic compressible fluids, Journal of Differential Equations., 190(2) (2003), 504--523.
       	\bibitem{13} M.T. Chen, B. Huang, J.W. Zhang, Blowup criterion for the three-dimensional equations of compressible viscous
       	micropolar fluids with vacuum, Nonlinear Anal., 79 (2013), 1-–11.
       	\bibitem{12} M.T. Chen, X.Y. Xu, J.W. Zhang, Global weak solutions of 3D compressible micropolar fluids with discontinuous
        initial data and vacuum, Commun. Math. Sci., 13 (2015), 225–-247.
       	\bibitem{8} B.Q. Dong, J.N. Li, J.H. Wu, Global well-posedness and large-time decay for the 2D micropolar equations, J. Differential Equations., 262 (2017), 3488–-3523.  
       	\bibitem{9} I. Dražić,  N. Mujaković, 3-D flow of a compressible viscous micropolar fluid with spherical symmetry: a global existence theorem, Boundary value problems., 1 (2015) 98.    
       	\bibitem{2} A.C. Eringen, Theory of micropolar fluids, J. Appl. Math. Mech., 16 (1966) 1–-18.
       	\bibitem{23} E. Feireisl,  A. Novotný, H. Petzeltová, On the existence of globally defined weak solutions to the Navier—Stokes equations, Journal of Mathematical Fluid Mechanics., 3(4) (2001), 358--392. 
       	\bibitem{18} X. Huang, J. Li, Z. P. Xin, Global well-posedness of classical solutions with large oscil-lations and vacuum to the three-dimensional isentropic compressible Navier-Stokes equations, Commun. Pure Appl. Math., 65 (2012), 549-–585.   
       	\bibitem{29} X. Huang, J. Li, Z. P. Xin, Serrin-type criterion for the three-dimensional viscous com-pressible flows, SIAM J. Math. Anal., 43 (2011), 1872–-1886.
       	\bibitem{3} G. Lukaszewicz, Micropolar Fluids: Theory and Applications, Birkhäuser, Boston., 1999.         
       \bibitem{22} P. L. Lions, Mathematical Topics in Fluid Mechanics, In Compressible models,Oxford University Press, New York., (vol. 2) 
       	(1998).           	
       	\bibitem{19} H. Li, X. Xu, J. Zhang, Global classical solutions to 3D compressible magnetohydrodynamic equations with large oscillations and vacuum, SIAM J. Math. Anal., 45 (2013), 1356--1387.   	
       	\bibitem{14} Q. Liu, P. Zhang, Optimal time decay of the compressible micropolar fluids, J. Differential Equations., 260 (10) (2016), 7634–-7661.    
       	\bibitem{6} N. Mujakovic, One-dimensional flow of a compressible viscous micropolar fluid: a global existence theorem, Glas.
       	Mat., 33 (53) (1998), 199–-208.       	
       	\bibitem{5} N. Mujakovic, One-dimensional flow of a compressible viscous micropolar fluid: a local existence theorem, Glas.
       	Mat., 33 (53) (1998), 71–-91.
       	\bibitem{7} N. Mujakovic, Non-homogeneous boundary value problem for one-dimensional compressible viscous micropolar
       	fluid model: a global existence theorem, Math. Inequal. Appl., 12 (2009), 651–-662.
       	\bibitem{10} N. Mujaković, I. Dražić, 3-D flow of a compressible viscous micropolar fluid with spherical symmetry: uniqueness of a generalized solution, Boundary value problems., 1 (2014), 226.       
       	\bibitem{21} A. Matsumura, T. Nishida, The initial value problem for the equations of motion of viscous and heatconductive gases, J Math Kyoto Univ., 20 (1980), 67-–104.
      	\bibitem{17} R. Wei, B. Guo, Y. Li, Global existence and optimal convergence rates of solutions for 3D compressible magneto-micropolar fluid equations, Journal of Differential Equations., 263(5) (2017), 2457--2480.	
      	\bibitem{15}  J. Yuan, Existence theorem and blow‐up criterion of the strong solutions to the magneto‐micropolar fluid equations, Mathematical Methods in the Applied Sciences., 31(9) (2008), 1113--1130.
       	\bibitem{27} A. A. Zlotnik, Uniform estimates and stabilization of symmetric solutions of a system of
       	quasilinear equations, Differ. Equ., 36 (2000), 701–-716.             
       	\end{thebibliography}
\end{document}